\documentclass[12pt]{amsart}
\usepackage{amsmath}
\usepackage{amsfonts}
\usepackage{tikz-cd}
\usepackage{amssymb}
\usepackage[all,cmtip]{xy}           
\usepackage{mathrsfs}
\usepackage{bbm}
\usepackage{bbding}
\usepackage{txfonts}
\usepackage[shortlabels]{enumitem}
\usepackage{ifpdf}
\ifpdf
\usepackage[colorlinks,final,backref=page,hyperindex]{hyperref}
\else
\usepackage[colorlinks,final,backref=page,hyperindex,hypertex]{hyperref}
\fi
\usepackage{amscd}
\usepackage{tikz-cd}
\usepackage[active]{srcltx}

\makeatletter

\newtheorem{theorem}{Theorem}[section]
\newtheorem{prop}[theorem]{Proposition}
\newtheorem{lemma}[theorem]{Lemma}
\newtheorem{coro}[theorem]{Corollary}
\newtheorem{prop-def}{Proposition-Definition}[section]

\theoremstyle{definition}
\newtheorem{defn}[theorem]{Definition}

\newtheorem{remark}[theorem]{Remark}
\newtheorem{exam}[theorem]{Example}

\newcommand{\nc}{\newcommand}

\nc{\delete}[1]{{}}
\nc{\mmargin}[1]{}

\nc{\mlabel}[1]{\label{#1}}  
\nc{\mcite}[1]{\cite{#1}}  
\nc{\mref}[1]{\ref{#1}}  
\nc{\mbibitem}[1]{\bibitem{#1}} 

\delete{
	\nc{\mlabel}[1]{\label{#1}  
		{\hfill \hspace{1cm}{\bf{{\ }\hfill(#1)}}}}
	\nc{\mcite}[1]{\cite{#1}{{\bf{{\ }(#1)}}}}  
	\nc{\mref}[1]{\ref{#1}{{\bf{{\ }(#1)}}}}  
	\nc{\mbibitem}[1]{\bibitem[\bf #1]{#1}} 
}

 \font\cyrs=wncyr7

\nc{\vep}{\varepsilon}
\nc{\bin}[2]{ (_{\stackrel{\scs{#1}}{\scs{#2}}})}  
\nc{\binc}[2]{(\!\! \begin{array}{c} \scs{#1}\\
		\scs{#2} \end{array}\!\!)}  
\nc{\bincc}[2]{  ( {\scs{#1} \atop
		\vspace{-1cm}\scs{#2}} )}  
\nc{\oline}[1]{\overline{#1}}
\nc{\mapm}[1]{\lfloor\!|{#1}|\!\rfloor}
\nc{\bs}{\bar{S}}
\nc{\la}{\longrightarrow}
\nc{\ot}{\otimes}
\nc{\rar}{\rightarrow}
\nc{\lon }{\,\rightarrow\,}
\nc{\dar}{\downarrow}
\nc{\dap}[1]{\downarrow \rlap{$\scriptstyle{#1}$}}
\nc{\defeq}{\stackrel{\rm def}{=}}
\nc{\dis}[1]{\displaystyle{#1}}
\nc{\dotcup}{\ \displaystyle{\bigcup^\bullet}\ }
\nc{\hcm}{\ \hat{,}\ }
\nc{\hts}{\hat{\otimes}}
\nc{\hcirc}{\hat{\circ}}
\nc{\lleft}{[}
\nc{\lright}{]}
\nc{\curlyl}{\left \{ \begin{array}{c} {} \\ {} \end{array}
	\right .  \!\!\!\!\!\!\!}
\nc{\curlyr}{ \!\!\!\!\!\!\!
	\left . \begin{array}{c} {} \\ {} \end{array}
	\right \} }
\nc{\longmid}{\left | \begin{array}{c} {} \\ {} \end{array}
	\right . \!\!\!\!\!\!\!}
\nc{\ora}[1]{\stackrel{#1}{\rar}}
\nc{\ola}[1]{\stackrel{#1}{\la}}
\nc{\scs}[1]{\scriptstyle{#1}} \nc{\mrm}[1]{{\rm #1}}
\nc{\dirlim}{\displaystyle{\lim_{\longrightarrow}}\,}
\nc{\invlim}{\displaystyle{\lim_{\longleftarrow}}\,}
\nc{\dislim}[1]{\displaystyle{\lim_{#1}}} \nc{\colim}{\mrm{colim}}
\nc{\mvp}{\vspace{0.3cm}} \nc{\tk}{^{(k)}} \nc{\tp}{^\prime}
\nc{\ttp}{^{\prime\prime}} \nc{\svp}{\vspace{2cm}}
\nc{\vp}{\vspace{8cm}}
\nc{\modg}[1]{\!<\!\!{#1}\!\!>}
\nc{\intg}[1]{F_C(#1)}
\nc{\lmodg}{\!<\!\!}
\nc{\rmodg}{\!\!>\!}
\nc{\cpi}{\widehat{\Pi}}
\nc{\ssha}{{\mbox{\cyrs X}}} 
\nc{\tsha}{{\mbox{\cyrt X}}}
\nc{\shpr}{\diamond}    
\nc{\labs}{\mid\!}
\nc{\rabs}{\!\mid}

\nc{\C}{{\mathrm{C}}}
\nc{\ad}{\mrm{ad}}
\nc{\ann}{\mrm{ann}}
\nc{\Aut}{\mrm{Aut}}
\nc{\DA}{{\mathsf{DL}_\lambda}}
\nc{\Alg}{{\mathrm{Lie}}}
\nc{\DO}{{\mathsf{DO}_\lambda}}
\nc{\bim}{\mbox{-}\mathsf{Rep}}
\nc{\md}{\mbox{-}\mathsf{rep}}
\nc{\br}{\mrm{bre}}
\nc{\can}{\mrm{can}}
\nc{\Cont}{\mrm{Cont}}
\nc{\rchar}{\mrm{char}}
\nc{\cok}{\mrm{coker}}
\nc{\de}{\mrm{dep}}
\nc{\dtf}{{R-{\rm tf}}}
\nc{\dtor}{{R-{\rm tor}}}

\nc{\Div}{{\mrm Div}}
\nc{\Diff}{\mrm{DL}}
\nc{\Diffl}{\mathsf{DL}_\lambda}
\nc{\diffo}{{\mathsf{DO}_\lambda}}
\nc{\Dif}{{\mathfrak{Dif}^\lambda}}
\nc{\Difinfty}{{\mathfrak{Dif}^\lambda_\infty}}
\nc{\alg}{\mathsf{Lie}}
\nc{\End}{\mrm{End}}
\nc{\Ext}{\mrm{Ext}}
\nc{\Fil}{\mrm{Fil}}
\nc{\Fr}{\mrm{Fr}}
\nc{\Frob}{\mrm{Frob}}
\nc{\Gal}{\mrm{Gal}}
\nc{\GL}{\mrm{GL}}
\nc{\Hom}{\mrm{Hom}}
\nc{\Hoch}{\mrm{Hoch}}
\nc{\hsr}{\mrm{H}}
\nc{\hpol}{\mrm{HP}}
\nc{\id}{\mrm{id}}
\nc{\im}{\mrm{im}}
\nc{\Id}{\mrm{Id}}
\nc{\ID}{\mrm{ID}}
\nc{\Irr}{\mrm{Irr}}
\nc{\incl}{\mrm{incl}}
\nc{\length}{\mrm{length}}
\nc{\NLSW}{\mrm{NLSW}}
\nc{\Lie}{\mrm{Lie}}
\nc{\mchar}{\rm char}
\nc{\mpart}{\mrm{part}}
\nc{\ql}{{\QQ_\ell}}
\nc{\qp}{{\QQ_p}}
\nc{\rank}{\mrm{rank}}
\nc{\rcot}{\mrm{cot}}
\nc{\rdef}{\mrm{def}}
\nc{\rdiv}{{\rm div}}
\nc{\rtf}{{\rm tf}}
\nc{\rtor}{{\rm tor}}
\nc{\res}{\mrm{res}}
\nc{\Sh}{{\mathrm{Sh}}}
\nc{\SL}{\mrm{SL}}
\nc{\Spec}{\mrm{Spec}}
\nc{\sgn}{{\mathrm{sgn}}}
\nc{\tor}{\mrm{tor}}
\nc{\Tr}{\mrm{Tr}}
\nc{\tr}{\mrm{tr}}
\nc{\wt}{\mrm{wt}}
\nc{\op}{\mrm{op}}
\nc{\rmH}{ {\mathrm{H}}}

\nc{\bfk}{{\bf k}}
\nc{\bfone}{{\bf 1}}
\nc{\bfzero}{{\bf 0}}
\nc{\detail}{\marginpar{\bf More detail}
	\noindent{\bf Need more detail!}
	\svp}
\nc{\gap}{\marginpar{\bf Incomplete}\noindent{\bf Incomplete!!}
	\svp}
\nc{\FMod}{\mathbf{FMod}}
\nc{\Int}{\mathbf{Int}}
\nc{\Mon}{\mathbf{Mon}}
\nc{\remarks}{\noindent{\bf Remarks: }}
\nc{\Rep}{\mathbf{Rep}}
\nc{\Rings}{\mathbf{Rings}}
\nc{\Sets}{\mathbf{Sets}}
\nc{\ob}{\mathsf{Ob}}

\nc{\BA}{{\mathbb A}}   \nc{\CC}{{\mathbb C}}
\nc{\DD}{{\mathbb D}}   \nc{\EE}{{\mathbb E}}
\nc{\FF}{{\mathbb F}}   \nc{\GG}{{\mathbb G}}
\nc{\HH}{{\mathbb H}}   \nc{\LL}{{\mathbb L}}
\nc{\NN}{{\mathbb N}}   \nc{\PP}{{\mathbb P}}
\nc{\QQ}{{\mathbb Q}}   \nc{\RR}{{\mathbb R}}
\nc{\TT}{{\mathbb T}}   \nc{\VV}{{\mathbb V}}
\nc{\ZZ}{{\mathbb Z}}   \nc{\TP}{\widetilde{P}}
\nc{\m}{{\mathbbm m}}

\nc{\cala}{{\mathcal A}}    \nc{\calc}{{\mathcal C}}
\nc{\cald}{\mathcal{D}}     \nc{\cale}{{\mathcal E}}
\nc{\calf}{{\mathcal F}}    \nc{\calg}{{\mathcal G}}
\nc{\calh}{{\mathcal H}}    \nc{\cali}{{\mathcal I}}
\nc{\call}{{\mathcal L}}    \nc{\calm}{{\mathcal M}}
\nc{\caln}{{\mathcal N}}    \nc{\calo}{{\mathcal O}}
\nc{\calp}{{\mathcal P}}    \nc{\calr}{{\mathcal R}}
\nc{\cals}{{\mathcal S}}    \nc{\calt}{{\Omega}}
\nc{\calv}{{\mathcal V}}    \nc{\calw}{{\mathcal W}}
\nc{\calx}{{\mathcal X}}

\nc{\fraka}{{\mathfrak a}}
\nc{\frakb}{\mathfrak{b}}
\nc{\frakg}{{\frak g}}
\nc{\frakl}{{\frak l}}
\nc{\fraks}{{\frak s}}
\nc{\frakB}{{\frak B}}
\nc{\frakm}{{\frak m}}
\nc{\frakM}{{\frak M}}
\nc{\frakp}{{\frak p}}
\nc{\frakW}{{\frak W}}
\nc{\frakX}{{\frak X}}
\nc{\frakS}{{\frak S}}
\nc{\frakA}{{\frak A}}
\nc{\frakx}{{\frakx}}
\nc{\frakC}{{\frak{C}}}

\nc{\lir}[1]{\textcolor{red}{\underline{Li:}#1 }}

\begin{document}

\title [Twisted Rota-Baxter  families on Lie-Yamaguti algebras ]{Twisted Rota-Baxter  families
on Lie-Yamaguti algebras and NS-Lie-Yamaguti family algebras}

\author{Wen Teng}\footnote{Corresponding author}
\address{ School of Mathematics and Statistics, Guizhou University of Finance and Economics,  Guiyang  550025, P. R. of China}
\email{tengwen@mail.gufe.edu.cn}

\date{\today}

\begin{abstract}
In this paper, we first introduce twisted  Rota-Baxter  families
on Lie-Yamaguti algebras indexed by a commutative semigroup $\Omega$. Then, we study  NS-Lie-Yamaguti family algebras as
the underlying structures of  twisted  Rota-Baxter families. Finally, we investigate the cohomology of a twisted  Rota-Baxter family.
This cohomology can also be seen as the cohomology of a certain $\Omega$-Lie-Yamaguti algebras with
coefficients in an appropriate representation. As  applications, we consider the
  deformations of twisted  Rota-Baxter  families from the cohomological points
of view.
\end{abstract}

\subjclass[2020]{
 17A40; 17B38; 17B56
}

\keywords{Lie-Yamaguti algebra; twisted   Rota-Baxter operator;  $\Omega$-Lie-Yamaguti algebra;  NS-Lie-Yamaguti family algebra;  cohomology; deformation}

\maketitle

\tableofcontents

\allowdisplaybreaks

\section{Introduction}
\def\theequation{\arabic{section}. \arabic{equation}}
\setcounter{equation} {0}

 A Lie-Yamaguti algebra traces its origins to Nomizu's seminal work on affine invariant connections in homogeneous spaces during the 1950s \cite{Nomizu1954}.
 In the later 1960s, Yamaguti introduced an algebraic structure, which he termed a general Lie triple system or Lie triple algebra \cite{Yamaguti57,Yamaguti67,Yamaguti60},
  as it serves as a generalization of both Lie algebra and Lie triple system.
  Kinyon and Weinstein initially referred to this object as a Lie-Yamaguti algebra during their study of Courant algebroids in the early 21st century \cite{Kinyon}.
  Since a Lie-Yamaguti algebra originates entirely from differential geometry and serves as a pivotal higher structure in mathematical physics,
  it has garnered substantial attention and has been extensively investigated in recent times.
  For instance, Benito and his collaborators extensively investigated irreducible Lie-Yamaguti algebras and their connections with orthogonal Lie algebras \cite{Benito11,Benito15,Benito05,Benito09}.
  The deformations and extensions of Lie-Yamaguti algebras were thoroughly examined in \cite{Goswami,Goswami24,Lin15,Zhang2015,Zhao2023}.
  Product structures and complex structures on Lie-Yamaguti algebras, as explored through the application of Nijenhuis operators in \cite{Sheng2021}.
  Rota-Baxter operators and symplectic structures on Lie-Yamaguti algebras have been studied in \cite{Sheng2022}.

The concept of Rota-Baxter operators originated from Baxter's work in 1960,
aimed at addressing issues in probability \cite{Baxter}, and was subsequently intensively investigated by Rota  \cite{Rota} in the context of combinatorics.
In the 1980s, the Rota-Baxter operator of weight 0 was introduced in relation to the classical Yang-Baxter equation for Lie algebras  \cite{Belavin}.
Subsequently, Kupershmidt \cite{Kupershmidt} introduced the concept of relative Rota-Baxter operators, also referred to as $\mathcal{O}$-operators.
This definition is intended to enhance the understanding of the classical Yang-Baxter equations and their related integrable systems.
A substantial body of scholarly research has been devoted to investigating the multifaceted characteristics of Rota-Baxter operators, encompassing both pure and applied mathematics \cite{Gubarev,Bai12,Ebrahimi-Fard07,Ebrahimi-Fard08}. For further information on Rota-Baxter operators, refer to \cite{Guo12}.

In \cite{Uchino}, Uchino introduced twisted  Rota-Baxter operators on an associative algebra. These operators are referred to as generalized Reynolds operators or twisted $\mathcal{O}$-operators, and they can be regarded as a non-commutative counterpart to twisted Poisson structures \cite{Severa}.
Later, the comprehensive examination of twisted Rota-Baxter operators on Lie algebras, their cohomology and deformation theory, as well as the associated NS-Lie algebras, is presented in \cite{Das2021}. Additionally, the thorough exploration of twisted Rota-Baxter operators on Lie-Yamaguti algebras, along with their cohomology and deformation theory, is detailed in \cite{Teng2023}.
For further reference, consult \cite{Mabrouk25}.

Algebraic structures frequently manifest in `family versions', where conventional operations are supplanted by a family of operations indexed by a semigroup $\Omega$.
The concept of family versions of algebraic structures initially emerged in the research of K. Ebrahimi-Fard, J. Gracia-Bondia, and F. Patras, particularly in the algebraic formulations of renormalization within quantum field theory \cite{Ebrahimi-Fard,Kreimer}. The concept of the Rota-Baxter family was also introduced in \cite{Guo2009}.
Recently, various family algebraic structures have been defined, as referenced in \cite{Foissy2021,Foissy2024,Zhang2019,Zhang2024}.
 In particular, Das \cite{Das22} introduced the concept of an $\mathcal{O}$-operator family on associative algebras, which serves as a generalization of the Rota-Baxter family. Furthermore, the cohomology and deformation theory of twisted $\mathcal{O}$-operator families on associative algebras, as well as NS-family algebras, have also been investigated.
Liu and Zheng studied the cohomology and deformation theory of  twisted   $\mathcal{O}$-operator families on    Leibniz algebras  in \cite{Liu}.
See also related research in \cite{Teng2025,Teng25,Teng24}.

Prompted by the previously mentioned findings, we delve into the study of twisted Rota-Baxter families in the context of Lie-Yamaguti algebras.
This paper is organized as follows.
In Section  \ref{sec:Preliminaries},
we recall some basic notions such as Lie-Yamaguti algebras, representations, and (2,3)-cocycle.
In Section  \ref{sec:Twisted Rota-Baxter}, we introduce twisted    Rota-Baxter families on  Lie-Yamaguti algebras, and give some examples.
In Section \ref{sec:Hom-NS-family algebras}, we introduce NS-Lie-Yamaguti   family algebras as the underlying structure of  twisted    Rota-Baxter families.
 We show that an NS-Lie-Yamaguti family algebra induces an ordinary NS-Lie-Yamaguti algebra
on the tensor product with the semigroup algebra.
 In Section \ref{sec:O-Cohomology},  we introduce $\Omega$-Lie-Yamaguti algebras and their cohomology.
 Finally, in Section \ref{sec:Cohomology}, we study cohomology and deformations of  twisted   Rota-Baxter families.

\section{Preliminaries   }\label{sec:Preliminaries}
\def\theequation{\arabic{section}.\arabic{equation}}
\setcounter{equation} {0}

In this section, we recall some basic notions such as Lie-Yamaguti algebras, representations, and (2,3)-cocycle, as presented in \cite{Kinyon} and \cite{Yamaguti67,Zhang2015}.
Throughout this paper, we work on an algebraically closed field $\mathbb{K}$ of characteristics different
from 2 and 3.

\begin{defn} (\cite{Kinyon})
A Lie-Yamaguti  algebra  is a  3-tuple $(L, [\cdot, \cdot], \{\cdot, \cdot, \cdot\})$ in which $L$ is a vector space together with  a bilinear bracket  $[\cdot, \cdot]: \wedge^2L\rightarrow L$ and  a  trilinear bracket $\{\cdot, \cdot, \cdot\}: \wedge^2L\otimes L\rightarrow L$ such that the following conditions hold
\begin{align}
&[[x, y], z]+[[y, z], x]+[[z, x], y]+\{x, y, z\}+\{z, x, y\}+\{y, z, x\}=0,\label{2.1}\\
& \{[x, y], z, a\}+ \{[y, z], x, a\}+ \{[z, x], y, a\}=0,\label{2.2}\\
&\{a, b, [x, y]\}=[\{a, b, x\}, y]+[x,\{a, b, y\}],\label{2.3}\\
& \{a, b, \{x, y, z\}\}=\{\{a, b, x\}, y, z\}+ \{x,  \{a, b, y\}, z\}+ \{x,  y, \{a, b, z\}\},\label{2.4}
\end{align}
for all $ x, y, z, a, b\in L$.
\end{defn}

A homomorphism between two Lie-Yamaguti  algebras  $(L, [\cdot, \cdot], \{\cdot, \cdot, \cdot\})$ and $(L', [\cdot, \cdot]',$ $ \{\cdot,\cdot, \cdot\}')$ is a linear map $\varphi: L\rightarrow L'$ satisfying
\begin{eqnarray*}
&&\varphi([x, y])=[\varphi(x), \varphi(y)]', ~~~\varphi(\{x, y, z\})=\{\varphi(x), \varphi(y),\varphi(z)\}',~~\forall ~x,y,z\in L.
\end{eqnarray*}

\begin{exam} \label{exam:Lie algebra}
 Let $(L, [\cdot,\cdot])$  be a Lie algebra.  Define a  ternary bracket on $L$ by
 $$\{x,y,z\}=[[x,y],z],\ \ \  \  \forall x,y,z\in L.$$
 Then  $(L, [\cdot, \cdot], \{\cdot, \cdot, \cdot\})$ becomes  a     Lie-Yamaguti  algebra.
\end{exam}

\begin{exam}  \label{exam:Leibniz algebra}
 Let $(L, \star)$  be a (left) Leibniz algebra.   Define a binary and ternary bracket on $L$ by
 $$[x,y]=x\star y-y\star x, \{x,y,z\}=-(x\star y)\star z, \forall x,y,z\in L.$$
 Then $(L, [\cdot, \cdot], \{\cdot, \cdot, \cdot\})$ is  a     Lie-Yamaguti  algebra.
\end{exam}

Yamaguti introduced the  representation of
Lie-Yamaguti   algebra in  \cite{Yamaguti67}, which is a natural generalization of that of a Lie algebra or a Lie triple system.

\begin{defn} (\cite{Yamaguti67} )
Let $(L, [\cdot, \cdot], \{\cdot, \cdot, \cdot\})$  be a  Lie-Yamaguti   algebra   and $V$ be a vector space. A representation of $(L, [\cdot, \cdot], \{\cdot, \cdot, \cdot\})$ on $V$ consists of a linear map  $\rho: L\rightarrow \mathrm{End}(V)$ and a bilinear map  $\theta: L\times L\rightarrow \mathrm{End}(V)$ such that $x,y,z,a,b\in L$,
\begin{align}
0=&\theta([x,y],z)-\theta(x,z)\rho(y)+\theta(y,z)\rho(x),\label{2.5}\\
0=&D(x,y)\rho(z)-\rho(z)D(x,y)-\rho(\{x,y,z\}),\label{2.6}\\
0=&\theta(x,[y,z])-\rho(y)\theta(x,z)+\rho(z)\theta(x,y),\label{2.7}\\
0=&D(a,b)\theta(x,y)-\theta(x,y)D(a,b)-\theta(\{a,b,x\},y)-\theta(x,\{a,b,y\}),\label{2.8}\\
0=&\theta(a,\{x,y,z\})-\theta(y,z)\theta(a,x)+\theta(x,z)\theta(a,y)-D(x,y)\theta(a,z),\label{2.9}
\end{align}
where  $D: L\times L$ to $\mathrm{End}(V)$ is given by
\begin{align}
D(x,y)=\theta(y,x)-\theta(x,y)-\rho([x,y])+\rho(x)\rho(y)-\rho(y)\rho(x).\label{2.10}
\end{align}
 In this case, we also call $V$ a $L$-module.
\end{defn}
It can be concluded from \eqref{2.5}-\eqref{2.10} that
\begin{align}
0=&D([x,y],z)+D([y,z],x)+D([z,x],y),\label{2.11}\\
0=&D(a,b)D(x,y)-D(x,y)D(a,b)-D(\{a,b,x\},y)-D(x,\{a,b,y\}),\label{2.12}\\
0=&\theta(\{x,y,z\},a)-\theta(x,a)\theta(z,y)+\theta(y,a)\theta(z,x)+\theta(z,a)D(x,y).\label{2.13}
\end{align}
\begin{exam}

Let $(L, [\cdot, \cdot], \{\cdot, \cdot, \cdot\})$  be a  Lie-Yamaguti   algebra.  We   define linear maps $\mathrm{ad}: L\rightarrow \mathrm{End}(L),\mathcal{R}:\otimes^2 L\rightarrow \mathrm{End}(L)$ by
\begin{eqnarray*}
& \mathrm{ad}(x)(z):=[x,z],\mathcal{R}(x,y)(z):=\{z,x,y\},
\end{eqnarray*}for all $x,y,z\in L$.  Then $(L;\mathrm{ad}, \mathcal{R})$ forms a representation of $L$ on itself,  called the adjoint representation,
where $\mathcal{L}:=D$
is given by
$$\mathcal{L}(x,y)(z)=\mathcal{R}(y,x)z-\mathcal{R}(x,y)z-\mathrm{ad}([x,y])z+\mathrm{ad}(x)\mathrm{ad}(y)z-\mathrm{ad}(y)\mathrm{ad}(x)z=\{x,y,z\}.$$
\end{exam}

Let $(V;\rho,\theta)$ be a representation of a Lie-Yamaguti algebra $(L, [\cdot, \cdot], \{\cdot, \cdot, \cdot\})$. Let us recall the (2,3)-cocycle on Lie-Yamaguti algebras in \cite{Yamaguti67,Zhang2015}.

In \cite{Zhang2015}, the authors employed the (2,3)-cocycle of Lie-Yamaguti algebra to classify and study  1-parameter infinitesimal deformations and abelian extensions of Lie-Yamaguti algebras.
A  (2,3)-cocycle on  a Lie-Yamaguti algebra $(L, [\cdot, \cdot], \{\cdot, \cdot, \cdot\})$ with coefficients in a representation $(V; \rho, \theta)$ is
a pair $\Gamma=(\Gamma_1,\Gamma_2)\in \mathrm{Hom}(\wedge^2 L,V)\times \mathrm{Hom}(\wedge^2 L\otimes L,V)$  that satisfies the following conditions:
\begin{small}
\begin{align}
0=&-\rho(x_1)\Gamma_1(y_1,z)-\rho(y_1)\Gamma_1(z,x_1)-\rho(z)\Gamma_1(x_1,y_1)+\Gamma_1([x_1,y_1],z)+\Gamma_1([y_1,z],x_1)+\label{2.14}\\
&\Gamma_1([z,x_1],y_1)+\Gamma_2(x_1,y_1,z)+\Gamma_2(y_1,z,x_1)+\Gamma_2(z,x_1,y_1),\nonumber\\
0=&\theta(x_1,y_2)\Gamma_1(y_1,x_2)+\theta(y_1,y_2)\Gamma_1(x_2,x_1)+\theta(x_2,y_2)\Gamma_1(x_1,y_1)+\Gamma_2([x_1,y_1],x_2,y_2)+\label{2.15}\\
 &\Gamma_2([y_1,x_2],x_1,y_2)+\Gamma_2([x_2,x_1],y_1,y_2),\nonumber\\
0=&-\rho(x_2)\Gamma_2(x_1,y_1,y_2)+\rho(y_2)\Gamma_2(x_1,y_1,x_2)+\Gamma_2(x_1,y_1,[x_2,y_2])+\label{2.16}\\
&D(x_1,y_1)\Gamma_1(x_2,y_2)-\Gamma_1(\{x_1,y_1,x_2\},y_2)-\Gamma_1(x_2,\{x_1,y_1,y_2\}),\nonumber\\
0=&-\theta(y_2,z)\Gamma_2(x_1,y_1,x_2)+\theta(x_2,z)\Gamma_2(x_1,y_1,y_2)+D(x_1,y_1)\Gamma_2(x_2,y_2,z)-D(x_2,y_2)\Gamma_2(x_1,y_1,z)-\label{2.17}\\
&\Gamma_2(\{x_1,y_1,x_2\},y_2,z)-\Gamma_2(x_2,\{x_1,y_1,y_2\}, z)+\Gamma_2(x_1, y_1, \{x_2,y_2,z\})-\Gamma_2(x_2,y_2,\{x_1,y_1,z\}),\nonumber
\end{align}
\end{small}
for all $x_1,y_1,x_2,y_2,z\in L.$
For the general cohomology theory of Lie-Yamaguti algebras, see \cite{Yamaguti67}.

%
%

\section{Twisted   Rota-Baxter families   }\label{sec:Twisted Rota-Baxter}
\def\theequation{\arabic{section}.\arabic{equation}}
\setcounter{equation} {0}

In this  section, let $\Omega$ be a fixed commutative semigroup,  we introduce the notion of   $\Gamma$-twisted   Rota-Baxter families on    Lie-Yamaguti  algebras.

\begin{defn}
Let $(V; \rho, \theta)$ be a representation of a Lie-Yamaguti algebra $(L, [\cdot, \cdot], \{\cdot, \cdot, \cdot\})$, and let $\Gamma=(\Gamma_1,\Gamma_2)$ be a (2,3)-cocycle.
A collection $\{T_\alpha:V\rightarrow L\}_{\alpha\in\Omega}$ of linear maps is said to be an
$\Gamma$-twisted     Rota-Baxter family    (on a representation $(V;\rho,\theta)$  over the Lie-Yamaguti  algebra $(L, [\cdot, \cdot], \{\cdot, \cdot, \cdot\})$) if $\{T_\alpha\}_{\alpha\in\Omega}$ satisfies
\begin{small}
\begin{align}
&[T_\alpha u, T_\beta v]=T_{\alpha\beta}\big(\rho(T_\alpha u) v-\rho(T_\beta v)u+\Gamma_1(T_\alpha u, T_\beta v)\big),\label{3.1}\\
&\{T_\alpha u, T_\beta v, T_\gamma w\}=T_{\alpha\beta\gamma}\big(D(T_\alpha u, T_\beta v)w-\theta(T_\alpha u, T_\gamma w)v+\theta(T_\beta v, T_\gamma w)u+\Gamma_2(T_\alpha u, T_\beta v,T_\gamma w)\big),\label{3.2}
\end{align}
\end{small}
for $u,v,w\in V$ and $\alpha,\beta,\gamma\in \Omega.$
\end{defn}

 Let $(L', [\cdot, \cdot]', \{\cdot, \cdot, \cdot\}')$  be another  Lie-Yamaguti   algebra and  $(V'; \rho', \theta')$ be a
representation of $L'$. Suppose that $\Gamma'=(\Gamma'_1,\Gamma'_2)\in  \mathrm{Hom}(\wedge^2 L',V')\times \mathrm{Hom}(\wedge^2 L'\otimes L',V')$ is a (2,3)-cocycle and $\{T'_\alpha:V'\rightarrow L'\}_{\alpha\in\Omega}$ is an $\Gamma'$-twisted    Rota-Baxter family.
A morphism of twisted  Rota-Baxter families from $\{T_\alpha:V\rightarrow L\}_{\alpha\in\Omega}$ to $\{T'_\alpha:V'\rightarrow L'\}_{\alpha\in\Omega}$ consists of a pair $(\zeta,\eta)$ of a
Lie-Yamaguti morphism
$\eta:(L, [\cdot, \cdot], \{\cdot, \cdot, \cdot\})\rightarrow (L', [\cdot, \cdot]', \{\cdot, \cdot, \cdot\}')$ and a linear map $\zeta:V\rightarrow V'$ satisfying
\begin{align}
&\eta\circ T_\alpha=T'_\alpha\circ\zeta,~~~~~\zeta\circ\Gamma_1=\Gamma_1'\circ(\eta\otimes\eta),~~~~\zeta\circ\Gamma_2=\Gamma_2'\circ(\eta\otimes\eta\otimes\eta), \label{3.3}\\
&\zeta(\rho(x)u)=\rho'(\eta(x))\zeta(u),~~~~\zeta(\theta(x,y)u)=\theta' (\eta(x),\eta(y))u,~~~~\forall x,y\in L, u\in V.\label{3.4}
\end{align}

\begin{remark}
 On the one hand, if a Lie-Yamaguti algebra $(L, [\cdot, \cdot], \{\cdot, \cdot, \cdot\})$  reduces
to a Lie triple system  $(L,\{\cdot,\cdot,\cdot\})$, that is $[\cdot,\cdot]=0$, we get
the notion of an $\Gamma_2$-twisted     Rota-Baxter family $\{T_\alpha:V\rightarrow L\}_{\alpha\in\Omega}$  on a Lie triple system  $(L,\{\cdot,\cdot,\cdot\})$ \cite{Teng2025}.
On the other hand, if a Lie-Yamaguti algebra
$(L, [\cdot, \cdot], \{\cdot, \cdot, \cdot\})$  reduces to a Lie algebra $(L,[\cdot,\cdot])$,
 that is $\{\cdot,\cdot,\cdot\}=0$, we get the notion of
an $\Gamma_1$-twisted     Rota-Baxter family $\{T_\alpha:V\rightarrow L\}_{\alpha\in\Omega}$ on a Lie algebra $(L,[\cdot,\cdot])$.
\end{remark}

Let $(L, [\cdot, \cdot], \{\cdot, \cdot, \cdot\})$  be a   Lie-Yamaguti   algebra.
A Reynolds family on $L$ consists of a collection $\{T_\alpha:L\rightarrow L\}_{\alpha\in\Omega}$  of linear maps satisfying:
\begin{align*}
[T_\alpha x, T_\beta y]=&T_{\alpha\beta}\big([T_\alpha x,y]+[x, T_\beta y]-[T_\alpha x,T_\beta y]\big),\\
\{T_\alpha x, T_\beta y, T_\gamma z\}=&T_{\alpha\beta\gamma}\big(\{T_\alpha x,T_\beta y,z\}+\{T_\alpha x, y, T_\gamma z\}+\{x,T_\beta y,T_\gamma z\}-\{T_\alpha x,T_\beta y,z\}\big),
\end{align*}
for all $x,y,z\in L.$

\begin{exam}  \label{exam:Reynolds family}
Let $(L, [\cdot, \cdot], \{\cdot, \cdot, \cdot\})$  be a   Lie-Yamaguti   algebra.
 Then $\Gamma_{\mathrm{ad}}=(-[\cdot, \cdot], -\{\cdot, \cdot, \cdot\})$ is a (2,3)-cocycle in the Yamaguti cohomology
of $L$ with coefficients in itself. With this notation, a Reynolds family $\{T_\alpha:L\rightarrow L\}_{\alpha\in\Omega}$ on $L$
is simply an $\Gamma_{\mathrm{ad}}$-twisted Rota-Baxter family.
\end{exam}

\begin{exam} \label{exam:generalized Reynolds operator} (see \cite{Teng2023})
A linear map $T:V\rightarrow L$ is said to a generalized Reynolds operator  if $T$ satisfies
\begin{align*}
[T u, T v]=&T\big(\rho(T u) v-\rho(T v)u+\Gamma_1(T u, T v)\big),\\
\{T u, T v, T w\}=&T\big(D(T u, T v)w-\theta(T u, T w)v+\theta(T v, T w)u+\Gamma_2(T u, T v,T w)\big),
\end{align*}
for $u,v,w\in V$.\\
Any generalized Reynolds operator on  a Lie-Yamaguti  algebra is  an $\Gamma$-twisted     Rota-Baxter family  indexed by the trivial semigroup reduced to one element.
\end{exam}

\begin{exam} \label{exam:2.5}
A relative Rota-Baxter  family  on a Lie-Yamaguti algebra $(L, [\cdot, \cdot], \{\cdot, \cdot, \cdot\})$  is a collection $\{T_\alpha:V\rightarrow L\}_{\alpha\in\Omega}$ of linear maps satisfying:
\begin{align*}
[T_\alpha u, T_\beta v]=&T_{\alpha\beta}\big(\rho(T_\alpha u) v-\rho(T_\beta v)u\big),\\
\{T_\alpha u, T_\beta v, T_\gamma w\}=&T_{\alpha\beta\gamma}\big(D(T_\alpha u, T_\beta v)w-\theta(T_\alpha u, T_\gamma w)v+\theta(T_\beta v, T_\gamma w)u\big),
\end{align*}
for $u,v,w\in V$ and $\alpha,\beta,\gamma\in \Omega.$ \\
Any relative Rota-Baxter  family  on a Lie-Yamaguti algebra is an   $\Gamma$-twisted     Rota-Baxter family  with
$\Gamma=0$.
\end{exam}

\begin{exam}
Let $\{T_\alpha:V\rightarrow L\}_{\alpha\in\Omega}$ be a  $\Gamma_1$-twisted     Rota-Baxter family
 on  a representation  $(V;\rho)$ over the Lie  algebra $(L,[\cdot,\cdot]$).
Then, $\{T_\alpha:V\rightarrow L\}_{\alpha\in\Omega}$  constitutes  an $(\Gamma_1,\Gamma_2)$-twisted    Rota-Baxter family on  a representation  $(V;\rho,\theta_\rho)$  over   the induced  Lie-Yamaguti algebra $(L, [\cdot, \cdot], \{\cdot, \cdot, \cdot\})$ given in Example \ref{exam:Lie algebra}. Here, $\theta_\rho(x,y)=\rho(y)\rho(x)$ and  $\Gamma_2(x,y,z)=\Gamma_1([x,y],z)-\rho(z)\Gamma_1(x,y)$, for all $x,y,z\in L.$
\end{exam}

\begin{exam} \label{exam: tensor product}
Let   $(L,[\cdot,\cdot], \{\cdot,\cdot,\cdot\})$  be  a  Lie-Yamaguti algebra. Then the tensor product $L\otimes \mathbb{K}\Omega$
is a Lie-Yamaguti algebra with the brackets $[\![\cdot,\cdot]\!]$, $\{\!\{\cdot,\cdot,\cdot\}\!\}$:
\begin{align*}
[\![a\otimes \alpha, b\otimes \beta]\!]=&[a, b]\otimes \alpha\beta,\\
\{\!\{a\otimes \alpha, b\otimes \beta, c\otimes \gamma\}\!\}=&\{a, b, c \}\otimes \alpha\beta\gamma,~\ ~ \forall a\otimes \alpha, b\otimes \beta, c\otimes \gamma\in L\otimes \mathbb{K}\Omega.
\end{align*}
Moreover, $L$ is a representation of $(L\otimes \mathbb{K}\Omega,[\![\cdot,\cdot]\!], \{\!\{\cdot,\cdot,\cdot\}\!\})$ with the  linear maps
\begin{align*}
&\hat{\rho}(a\otimes \alpha)b =[a, b],~~\hat{\theta}(a\otimes \alpha,b\otimes \beta)c=\{c,a,b\}.
\end{align*}
Then we can prove that the map $\hat{\Gamma}=(\hat{\Gamma}_1,\hat{\Gamma}_2)\in\mathrm{Hom}\big(\wedge^2 (L\otimes \mathbb{K}\Omega),L\big)\times \mathrm{Hom}\big(\wedge^2 (L\otimes \mathbb{K}\Omega)\otimes (L\otimes \mathbb{K}\Omega),L\big)$ given by
$\big(\hat{\Gamma}_1(a\otimes \alpha,b\otimes \beta),\hat{\Gamma}_2(a\otimes \alpha,b\otimes \beta,c\otimes \gamma)\big)=(-[a, b],-\{a,b,c\})$ is a (2,3)-cocycle in the cohomology of $L\otimes \mathbb{K}\Omega$ with coefficients in
$L$. And the collection $\{\mathrm{id}_\alpha:L\rightarrow  L\otimes \mathbb{K}\Omega\}_{\alpha\in\Omega}$ defined by $\mathrm{id}_\alpha(a)=a\otimes\alpha$ for
$a\in L$ is an $\hat{\Gamma}$-twisted  Rota-Baxter family.
\end{exam}

\begin{exam}  \label{exam:Nijenhuis family}
Let   $(L,[\cdot,\cdot], \{\cdot,\cdot,\cdot\})$  be  a  Lie-Yamaguti algebra,  and $\{N_\alpha:L\rightarrow L\}_{\alpha\in\Omega}$ a Nijenhuis family on it, i.e.
\begin{align*}
[N_\alpha x,  N_\beta y]=&N_{\alpha\beta}\big( [N_\alpha x, y]+[x,  N_\beta y]- N_{\alpha\beta}([x, y])\big),\\
\{N_\alpha x,  N_\beta y,  N_\gamma z \}=&N_{\alpha\beta\gamma}\big(\{x, N_\beta y, N_\gamma z\}+\{N_\alpha x,  y, N_\gamma z\}+\{N_\alpha x,  N_\beta y, z\}\big)- \\
&N_{\alpha\beta\gamma}^2\big(\{N_\alpha x, y, z\}+\{x,  N_\beta y,z\}+\{x,  y, N_\beta z\}\big)+N_{\alpha\beta\gamma}^3\{x, y,  z\},~~  \forall x,y,z\in L, \alpha,\beta,\gamma\in\Omega.
\end{align*}
In this case, $L\otimes \mathbb{K}\Omega$ carries a new  Lie-Yamaguti algebra structure with the brackets $[\cdot,\cdot]_N$, $\{\cdot,\cdot,\cdot\}_N$:
\begin{align*}
[x\otimes \alpha, y\otimes \beta]_N=&\big( [N_\alpha x, y]+[x,  N_\beta y]- N_{\alpha\beta}([x, y])\big)\otimes \alpha\beta,\\
\{x\otimes \alpha, y\otimes \beta, z\otimes \gamma\}_N=&\Big(\big(\{x, N_\beta y, N_\gamma z\}+\{N_\alpha x,  y, N_\gamma z\}+\{N_\alpha x,  N_\beta y, z\}\big)- \\
&N_{\alpha\beta\gamma}\big(\{N_\alpha x, y, z\}+\{x,  N_\beta y,z\}+\{x,  y, N_\beta z\}\big)+N_{\alpha\beta\gamma}^2\{x, y,  z\}\Big)\otimes \alpha\beta\gamma,
\end{align*}
 for all $x,y,z\in L$ and $\alpha, \beta, \gamma\in \Omega.$  Moreover, $L$ is a representation of $(L\otimes \mathbb{K}\Omega,[\cdot,\cdot]_N, \{\cdot,\cdot,\cdot\}_N)$ with the linear maps
\begin{align*}
&\rho_N(x\otimes \alpha) y =[N_\alpha x, y],~~ \theta_N(x\otimes \alpha, y\otimes \beta)z=\{z, N_\alpha x, N_\beta y \}.
\end{align*}
It is easy to prove that the map  $\Gamma_N=(\Gamma_1^N, \Gamma_2^N)\in\mathrm{Hom}\big(\wedge^2 (L\otimes \mathbb{K}\Omega),L\big)\times \mathrm{Hom}\big(\wedge^2 (L\otimes \mathbb{K}\Omega)\otimes (L\otimes \mathbb{K}\Omega),L\big)$ given by
$\big(\Gamma_1^N(x\otimes \alpha,y\otimes \beta), \Gamma_2^N(x\otimes \alpha,y\otimes \beta,z\otimes \gamma)\big)=\Big(-N_{\alpha\beta}[x, y],-N_{\alpha\beta\gamma}\big(\{N_\alpha x, y, z\}+\{x,  N_\beta y,z\}+\{x,  y, N_\beta z\}-N_{\alpha\beta\gamma}\{x, y,  z\} \big)\Big)$ is a (2,3)-cocycle in the cohomology of $L\otimes \mathbb{K}\Omega$ with coefficients in
$L$. And the collection $\{\mathrm{id}_\alpha:L\rightarrow  L\otimes \mathbb{K}\Omega\}_{\alpha\in\Omega}$ defined by $\mathrm{id}_\alpha(x)=x\otimes\alpha$ for
$x\in L$ is an $\Gamma_N$-twisted  Rota-Baxter family.
\end{exam}

Let $(L, [\cdot, \cdot], \{\cdot, \cdot, \cdot\})$  be a   Lie-Yamaguti   algebra and  $(V; \rho, \theta)$ be a representation of $L$.
Consider the Lie-Yamaguti algebra $(L\otimes \mathbb{K}\Omega,[\![\cdot,\cdot]\!], \{\!\{\cdot,\cdot,\cdot\}\!\})$  given in Example \ref{exam: tensor product}, we can prove that
$V\otimes \mathbb{K}\Omega$  is a  representation  of $(L\otimes \mathbb{K}\Omega,[\![\cdot,\cdot]\!], \{\!\{\cdot,\cdot,\cdot\}\!\})$ with the  bilinear maps
\begin{align*}
&\bar{\rho} (x\otimes \alpha)(u\otimes \beta) =\rho(x)u \otimes\alpha\beta,~~\bar{\theta}(x\otimes \alpha, y\otimes \beta)(u\otimes \gamma)=\theta(x,y)u\otimes  \alpha\beta\gamma.
\end{align*}
Let  $\Gamma=(\Gamma_1,\Gamma_2)$ be a (2,3)-cocycle in the cohomology of $(L, [\cdot, \cdot], \{\cdot, \cdot, \cdot\})$  with
coefficients in $(V; \rho, \theta)$.
Then $\hat{\Gamma}=(\bar{\Gamma}_1,\bar{\Gamma}_2)\in\mathrm{Hom}\big(\wedge^2 (L\otimes \mathbb{K}\Omega),V\otimes \mathbb{K}\Omega\big)\times \mathrm{Hom}\big(\wedge^2 (L\otimes \mathbb{K}\Omega)\otimes (L\otimes \mathbb{K}\Omega),V\otimes \mathbb{K}\Omega\big)$ given by
$\big(\bar{\Gamma}_1(a\otimes \alpha,b\otimes \beta),\bar{\Gamma}_2(a\otimes \alpha,b\otimes \beta,c\otimes \gamma)\big)=\big(\Gamma_1(a,b),\Gamma_2(a,b,c)\big)\otimes \alpha\beta\gamma$ is a (2,3)-cocycle in the cohomology of $L\otimes \mathbb{K}\Omega$ with coefficients in
$V\otimes \mathbb{K}\Omega$. Under
this assumption, we have the following conclusion.

 \begin{prop} \label{prop:2.9}
Let $\{T_\alpha:V\rightarrow L\}_{\alpha\in\Omega}$  be an $\Gamma$-twisted   Rota-Baxter family.
Then the map
\begin{align*}
& \bar{T}:V\otimes \mathbb{K}\Omega\rightarrow L\otimes \mathbb{K}\Omega,~~ \bar{T}(u\otimes\alpha)= T_\alpha u\otimes\alpha
\end{align*}
is a generalized Reynolds operator  on $V\otimes \mathbb{K}\Omega$ over the Lie-Yamaguti algebra $(L\otimes \mathbb{K}\Omega,[\![\cdot,\cdot]\!], \{\!\{\cdot,\cdot,\cdot\}\!\})$.
 \end{prop}

  \begin{proof}
  For all $u,v,w\in V$ and $\alpha,\beta,\gamma\in \Omega,$ we obtain
\begin{align*}
&[\![\bar{T} ( u\otimes \alpha),  \bar{T}  (v\otimes \beta)]\!]=[\![T_\alpha u\otimes\alpha, T_\beta v\otimes \beta]\!]\\
&=[T_\alpha u,  T_\beta v]\otimes \alpha\beta=T_{\alpha\beta}(\rho(T_\alpha u)v-\rho(T_\beta v)u+\Gamma_1(T_\alpha u, T_\beta v))\otimes \alpha\beta\\
&=\bar{T} \big((\rho(T_\alpha u)v-\rho(T_\beta v)u+\Gamma_1(T_\alpha u, T_\beta v))\otimes \alpha\beta\big)\\
&=\bar{T} \big(\rho(T_\alpha u) v \otimes \alpha\beta-\rho(T_\beta v)u \otimes\alpha\beta+\Gamma_1(T_\alpha u, T_\beta v ) \otimes \alpha\beta\big)\\
&=\bar{T}\Big(\bar{\rho}\big(\bar{T}(u\otimes \alpha)\big)(v\otimes \beta)-\bar{\rho}\big(\bar{T} (v\otimes \beta)\big)( u\otimes \alpha)+\bar{\Gamma}_1\big(\bar{T}(u\otimes \alpha), \bar{T}(v\otimes \beta)\big)\Big)
\end{align*}
and
\begin{align*}
&\{\!\{\bar{T} ( u\otimes \alpha),  \bar{T}  (v\otimes \beta),  \bar{T}  (w\otimes \gamma)\}\!\}=\{\!\{T_\alpha u\otimes\alpha, T_\beta v\otimes \beta, T_\gamma\otimes\gamma \}\!\}\\
=&\{T_\alpha u,  T_\beta v, T_\gamma w\}\otimes \alpha\beta\gamma\\
=&T_{\alpha\beta\gamma}\big(D(T_\alpha u, T_\beta v)w-\theta(T_\alpha u, T_\gamma w)v+\theta(T_\beta v, T_\gamma w)u+\Gamma_2(T_\alpha u, T_\beta v,T_\gamma w)\big)\otimes \alpha\beta\gamma\\
=&\bar{T} \big((D(T_\alpha u, T_\beta v)w-\theta(T_\alpha u, T_\gamma w)v+\theta(T_\beta v, T_\gamma w)u+\Gamma_2(T_\alpha u, T_\beta v,T_\gamma w))\otimes \alpha\beta\gamma\big)\\
=&\bar{T} \big((D(T_\alpha u, T_\beta v)w\otimes \alpha\beta\gamma-\theta(T_\alpha u, T_\gamma w)v\otimes \alpha\beta\gamma+\theta(T_\beta v, T_\gamma w)u\otimes \alpha\beta\gamma+\Gamma_2(T_\alpha u, T_\beta v,T_\gamma w) \otimes \alpha\beta\gamma\big)\\
=&\bar{T} \Big(D\big(\bar{T} (u\otimes\alpha), \bar{T} (v\otimes\beta)\big)(w\otimes\gamma)-\theta\big(\bar{T}(u\otimes\alpha), \bar{T}(w\otimes\gamma)\big)(v\otimes\beta)+\theta\big(\bar{T}(v\otimes \beta), \bar{T}(w\otimes\gamma)\big)(u\otimes\alpha)+\\
&\bar{\Gamma}_2\big(\bar{T}(u\otimes \alpha), \bar{T}(v\otimes \beta),\bar{T}(w\otimes \gamma)\big)\Big).
\end{align*}
Therefore, $\bar{T}$ is a generalized Reynolds operator.
   \end{proof}

Next, we characterize an $\Gamma$-twisted   Rota-Baxter family   by its graph.
 Let $(L, [\cdot, \cdot], \{\cdot, \cdot, \cdot\})$  be a   Lie-Yamaguti   algebra and  $(V; \rho, \theta)$ be a representation of it.
  Suppose that $\Gamma$
is a (2,3)-cocycle, then the direct sum $L\oplus V$ carries a  Lie-Yamaguti   algebra structure with
the  brackets $[\cdot,\cdot]_\Gamma$  and $\{\cdot,\cdot,\cdot\}$ given by
\begin{align*}
[(x,u),(y,v)]_\Gamma=&\big([x, y], \rho(x)v- \rho(y)u+\Gamma_1(x,y)\big),\\
\{(x,u),(y,v),(z,w)\}_\Gamma=&\big(\{x, y, z\}, D(x,y)w- \theta(x,z)v+\theta(y,z)u+\Gamma_2(x,y,z)\big),
\end{align*}
for  all $x,y,z\in L, u,v,w\in V.$  This is called the $\Gamma$-twisted semidirect product  Lie-Yamaguti    algebra, which is denoted
by $L\ltimes_\Gamma V$.

Let  $(L, [\cdot, \cdot], \{\cdot, \cdot, \cdot\})$  be a   Lie-Yamaguti   algebra. A collection $\{M_\alpha\}_{\alpha\in \Omega}$ of subspaces of $L$ is said to
be a  subalgebra family of $L$ if   $[M_\alpha, M_\beta]\subseteq M_{\alpha\beta}$ and $\{M_\alpha, M_\beta, M_\gamma\}\subseteq M_{\alpha\beta\gamma}$ for $\alpha,\beta,\gamma\in \Omega$.

 \begin{prop}
A collection $\{T_\alpha:V\rightarrow P\}_{\alpha\in\Omega}$  of linear maps is an $\Gamma$-twisted    Rota-Baxter family    if and only if the collection of graphs
$\{Gr(T_\alpha)\}_{\alpha\in\Omega}=\{(T_\alpha u,u)~|~u\in V\}_{\alpha\in\Omega}$ is a subalgebra family of the $\Gamma$-twisted semidirect product  Lie-Yamaguti    algebra $L\ltimes_\Gamma V$.
 \end{prop}

 \begin{proof}
For any $(T_\alpha u,u)\in Gr(T_\alpha),(T_\beta v,v)\in Gr(T_\beta), (T_\gamma w,w)\in Gr(T_\gamma)$, we have
\begin{align*}
[(T_\alpha u,u),(T_\beta v,v)]_\Gamma=&([T_\alpha u, T_\beta v], \rho(T_\alpha u)v- \rho(T_\beta v)u+\Gamma_1(T_\alpha u,T_\beta v))
\end{align*}
and
\begin{align*}
&\{(T_\alpha u,u),(T_\beta v,v),(T_\gamma w,w)\}_\Gamma\\
=&\big(\{T_\alpha u, T_\beta v, T_\gamma w\}, D(T_\alpha u,T_\beta v)w- \theta(T_\alpha u,T_\gamma w)v+\theta(T_\beta v,T_\gamma w)u+\Gamma_2(T_\alpha u,T_\beta v,T_\gamma w)\big).
\end{align*}
Then  the collection of graphs
$\{Gr(T_\alpha)\}_{\alpha\in\Omega}$ is a subalgebra family of   $L\ltimes_\Gamma V$
if and only if
\begin{align*}
&[T_\alpha u, T_\beta v]=T_{\alpha\beta}\big(\rho(T_\alpha u)v- \rho(T_\beta v)u+\Gamma_1(T_\alpha u,T_\beta v)\big) ~\ ~\ ~~\text{and}~\\
&\{T_\alpha u, T_\beta v, T_\gamma w\}=T_{\alpha\beta\gamma}\big(D(T_\alpha u,T_\beta v)w- \theta(T_\alpha u,T_\gamma w)v+\theta(T_\beta v,T_\gamma w)u+\Gamma_2(T_\alpha u,T_\beta v,T_\gamma w)\big),
\end{align*}
that is, the collection $\{T_\alpha:V\rightarrow L\}_{\alpha\in\Omega}$   is an $\Gamma$-twisted   Rota-Baxter family.
 \end{proof}


\section{NS-Lie-Yamaguti family algebras} \label{sec:Hom-NS-family algebras}
\def\theequation{\arabic{section}.\arabic{equation}}
\setcounter{equation} {0}

In this section, we introduce NS-Lie-Yamaguti family algebras as the family analogue of NS-Lie-Yamaguti algebras.   We
observed that an   NS-Lie-Yamaguti family algebra induces an ordinary NS-Lie-Yamaguti  algebra. Finally,  we  show  that
a twisted  Rota-Baxter family  on  a  Lie-Yamaguti  algebra induces an NS-Lie-Yamaguti family algebra.

\begin{defn}
An NS-Lie-Yamaguti algebra is a quintuple $(L,\bullet,\vee, [\cdot, \cdot, \cdot],\{\cdot, \cdot, \cdot\})$ consisting of  a vector space $L$,   two bilinear operations $\bullet$ and $\vee:L\times L\rightarrow L$,
 where $\vee$ is skew-symmetric, and
two multilinear ternary operations $[\cdot, \cdot, \cdot]$ and $\{\cdot, \cdot, \cdot\}:L\times L\times L\rightarrow L$, with the property that  $ [x,y,z]=-[y,x,z]$ for any $x,y,z\in L$.
 Additionally, these operations satisfy the following conditions for all elements $x,y,z,a,b\in L$:
 \begin{align}
0=&\{a,[x,y]^\ast, z\}-\{y\bullet a, x,z\}+\{x\bullet a, y,z\},\label{4.1}\\
0=&\{x,y,z\bullet a\}^\ast-z\bullet\{x,y,a\}^\ast-[\![x,y,z]\!]\bullet a,\label{4.2}\\
0=&\{a,x,[y,z]^\ast\}-y\bullet\{a,x,z\}+z\bullet\{a,x,y\},\label{4.3}\\
0=& \{x,y,\{b, z, a\}\}^\ast-\{\{x,y,b\}^\ast, z, a\}-\{b, [\![x,y,z]\!], a\}-
\{b, z, [\![x,y,a]\!]\},\label{4.4}\\
0=&\{b, x, [\![y,z,a]\!]\}-\{\{b, x, y\},z, a\}+\{\{b, x, z\}, y, a\}
- \{y,z,\{b, x, a\}\}^\ast,\label{4.5}\\
0=&[x, y]^\ast  \vee  z+[y, z]^\ast  \vee  x+[z, x]^\ast \vee  y-z\bullet  (x \vee  y)-x\bullet  (y \vee  z)-y\bullet  (z \vee  x)+\label{4.6}\\
&[x,y,z] +[y,z,x] +[z,x,y] , \nonumber\\
 0=&\{x \vee  y,z,a\} +\{y \vee  z,x,a\} +\{z \vee  x,y,a\} +[[x, y]^\ast,z,a]+[[y, z]^\ast,x,a]+[[z, x]^\ast,y,a],\label{4.7}\\
0=&-z\bullet [x, y, a] +a\bullet  [x, y, z] +[x, y, [z, a]^\ast ] +\{x,y,z\vee  a\}^\ast -[\![x,y,z]\!] \vee  a-\label{4.8}\\
&z\vee  [\![x,y,a]\!], \nonumber\\
0=& -\{[x, y, z], a, b\}+\{[x, y, a], z, b\}+\{x,y,[z, a, b]\}^\ast-\{z,a,[x, y, b]\}^\ast-\label{4.9} \\
&[[\![x,y,z]\!], a, b]-[z, [\![x,y,a]\!], b]+[x, y, [\![z,a,b]\!]]-[z, a, [\![x,y,b]\!]],\nonumber
\end{align}
where
 \begin{align}
 [x, y]^\ast=&x\bullet y-y\bullet x+x \vee y,\label{4.10}\\
 \{x,y,z\}^\ast=&\{z,y,x\}-\{z,x,y\}+x\bullet (y\bullet z)-y\bullet (x\bullet z)-[x,y]^\ast\bullet z,\label{4.11}\\
 [\![x,y,z]\!]=&\{x,y,z\}^\ast+\{x,y,z\}-\{y,x,z\}+[x,y,z].\label{4.12}
 \end{align}
\end{defn}

\begin{remark}
Let $(L,\bullet,\vee, [\cdot, \cdot, \cdot],\{\cdot, \cdot, \cdot\})$ be an  NS-Lie-Yamaguti algebra.
Then the following equalities are satisfied:
 \begin{align}
0=&  \{[x, y]^\ast,z,a\}^\ast+\{[y, z]^\ast,x,a\}^\ast+\{[z, x]^\ast,y,a\}^\ast,\label{4.13}\\
0=&\{x,y,\{z,a,b\}^\ast\}^\ast-\{z,a,\{x,y,b\}^\ast\}^\ast-\{[\![x,y,z]\!],a,b\}^\ast-
\{z,[\![x,y,a]\!],b\}^\ast, \label{4.14}\\
0=&\{b, [\![x,y,z]\!], a\}-\{\{b, z, y\}, x, a\}+\{\{b, z, x\}, y, a\}+\{\{x, y, b\}^\ast, z, a\}.\label{4.15}
 \end{align}
\end{remark}

\begin{remark}
 If the trilinear operation $\{\cdot, \cdot, \cdot\}$ and the bilinear operation $\bullet$ in the above definition are trivial, then $(L, \vee, [\cdot, \cdot, \cdot])$
becomes a Lie-Yamaguti algebra.
Additional, if $\vee$ and $[\cdot, \cdot, \cdot]$ are trivial, then $(L,\bullet, \{\cdot, \cdot, \cdot\})$ is a pre-Lie-Yamaguti algebra \cite{Sheng2022}.
Therefore, NS-Lie-Yamaguti algebras are a generalization of both Lie-Yamaguti algebras
and pre-Lie-Yamaguti algebras.
\end{remark}

 The family version of  NS-Lie-Yamaguti algebra  is given by the following.

\begin{defn}
  An   NS-Lie-Yamaguti family algebra is  a  linear space $L$ equipped with
a family of bilinear operations  $\{\bullet_{\alpha},\vee_{\alpha,\beta}:L\otimes L\rightarrow L\}_{\alpha,\beta\in \Omega}$, where
$x\vee_{\alpha,\beta}y=-y\vee_{\beta,\alpha}x$,
and a family of multilinear  operations  $\big\{\{\cdot, \cdot, \cdot\}_{\beta,\gamma},[\cdot, \cdot, \cdot]_{\alpha,\beta,\gamma}:L\otimes L\otimes L\rightarrow L\big\}_{\alpha,\beta,\gamma\in \Omega}$ ,where
$[x,y,z]_{\alpha,\beta,\gamma}=-[y, x, z]_{\beta,\alpha,\gamma}$,
 that satisfy the following  equations:
 \begin{small}
 \begin{align}
0=&\{a,[x,y]^\ast_{\alpha,\beta}, z\}_{\alpha\beta,\gamma}-\{y\bullet_\beta a, x,z\}_{\alpha,\gamma}+\{x\bullet_\alpha a, y,z\}_{\beta,\gamma},\label{4.16}\\
0=&\{x,y,z\bullet_\gamma a\}^\ast_{\alpha,\beta}-z\bullet_\gamma\{x,y,a\}^\ast_{\alpha,\beta}-[\![x,y,z]\!]_{\alpha,\beta,\gamma}\bullet_{\alpha\beta\gamma} a,\label{4.17}\\
0=&\{a,x,[y,z]^\ast_{\beta,\gamma}\}_{\alpha,\beta\gamma}-y\bullet_\beta\{a,x,z\}_{\alpha,\gamma}+z\bullet_\gamma\{a,x,y\}_{\alpha,\beta},\label{4.18}\\
0=& \{x,y,\{b, z, a\}_{\gamma,\varsigma}\}^\ast_{\alpha,\beta}-\{\{x,y,b\}^\ast_{\alpha,\beta}, z, a\}_{\gamma,\varsigma}-\{b, [\![x,y,z]\!]_{\alpha,\beta,\gamma}, a\}_{\alpha\beta\gamma,\varsigma}-
\{b, z, [\![x,y,a]\!]_{\alpha,\beta,\varsigma}\}_{\gamma,\alpha\beta\varsigma},\label{4.19}\\
0=&\{b, x, [\![y,z,a]\!]_{\beta,\gamma,\varsigma}\}_{\alpha,\beta\gamma\varsigma}-\{\{b, x, y\}_{\alpha,\beta},z, a\}_{\gamma,\varsigma}+\{\{b, x, z\}_{\alpha,\gamma}, y, a\}_{\beta,\varsigma}
- \{y,z,\{b, x, a\}_{\alpha,\varsigma}\}^\ast_{\beta,\gamma},\label{4.20}\\
0=&[x, y]^\ast_{\alpha,\beta} \vee_{\alpha\beta,\gamma} z+[y, z]^\ast_{\beta,\gamma} \vee_{\beta\gamma,\alpha} x+[z, x]^\ast_{\gamma,\alpha} \vee_{\gamma\alpha,\beta} y-z\bullet_\gamma (x \vee_{\alpha,\beta} y)-x\bullet_\alpha (y \vee_{\beta,\gamma} z)-\label{4.21}\\
&y\bullet_\beta (z \vee_{\gamma,\alpha} x)+[x,y,z]_{\alpha,\beta,\gamma}+[y,z,x]_{\beta,\gamma,\alpha}+[z,x,y]_{\gamma,\alpha,\beta}, \nonumber\\
 0=&\{x \vee_{\alpha,\beta} y,z,a\}_{\gamma,\varsigma}+\{y \vee_{\beta,\gamma} z,x,a\}_{\alpha,\varsigma}+\{z \vee_{\gamma,\alpha} x,y,a\}_{\beta,\varsigma}+[[x, y]^\ast_{\alpha,\beta},z,a]_{\gamma,\varsigma}+\label{4.22}\\
&[[y, z]^\ast_{\beta,\gamma},x,a]_{\alpha,\varsigma}+[[z, x]^\ast_{\gamma,\alpha},y,a]_{\beta,\varsigma},\nonumber\\
0=&-z\bullet_\gamma [x, y, a]_{\alpha,\beta,\varsigma}+a\bullet_\varsigma [x, y, z]_{\alpha,\beta,\gamma}+[x, y, [z, a]^\ast_{\gamma,\varsigma}]_{\alpha,\beta,\gamma\varsigma}+\label{4.23}\\
&\{x,y,z\vee_{\gamma,\varsigma} a\}^\ast_{\alpha,\beta}-[\![x,y,z]\!]_{\alpha,\beta,\gamma}\vee_{\alpha\beta\gamma,\varsigma} a-z\vee_{\gamma,\alpha\beta\varsigma} [\![x,y,a]\!]_{\alpha,\beta,\varsigma},\nonumber\\
0=& -\{[x, y, z]_{\alpha,\beta,\gamma}, a, b\}_{\varsigma,\tau}+\{[x, y, a]_{\alpha,\beta,\varsigma}, z, b\}_{\gamma,\tau}+\{x,y,[z, a, b]_{\gamma,\varsigma,\tau}\}^\ast_{\alpha,\beta}-\label{4.24}\\
&\{z,a,[x, y, b]_{\alpha,\beta,\tau}\}^\ast_{\gamma,\varsigma}-[[\![x,y,z]\!]_{\alpha,\beta,\gamma}, a, b]_{\alpha\beta\gamma,\varsigma,\tau}-[z, [\![x,y,a]\!]_{\alpha,\beta,\varsigma}, b]_{\gamma,\alpha\beta\varsigma,\tau}+\nonumber\\
&[x, y, [\![z,a,b]\!]_{\gamma,\varsigma,\tau}]_{\alpha,\beta,\gamma\varsigma\tau}-[z, a, [\![x,y,b]\!]_{\alpha,\beta,\tau}]_{\gamma,\varsigma,\alpha\beta\tau},\nonumber
\end{align}
\end{small}
where $\big\{[\cdot, \cdot]^\ast_{\alpha,\beta},\{\cdot, \cdot, \cdot\}^\ast_{\alpha,\beta}, [\![\cdot, \cdot, \cdot]\!]_{\alpha,\beta,\gamma}\big\}_{\alpha,\beta,\gamma\in \Omega}$ are defined by
 \begin{align}
 [x, y]^\ast_{\alpha,\beta}=&x\bullet_\alpha y-y\bullet_\beta x+x \vee_{\alpha,\beta} y,\label{4.25}\\
 \{x,y,z\}^\ast_{\alpha,\beta}=&\{z,y,x\}_{\beta,\alpha}-\{z,x,y\}_{\alpha,\beta}+x\bullet_\alpha (y\bullet_\beta z)-y\bullet_\beta (x\bullet_\alpha z)-[x,y]^\ast_{\alpha,\beta}\bullet_{\alpha\beta} z,\label{4.26}\\
 [\![x,y,z]\!]_{\alpha,\beta,\gamma}=&\{x,y,z\}^\ast_{\alpha,\beta}+\{x,y,z\}_{\beta,\gamma}-\{y,x,z\}_{\alpha,\gamma}+[x,y,z]_{\alpha,\beta,\gamma},\label{4.27}
 \end{align}
  for all  $x,y,z,a, b\in L$  and  $\alpha,\beta,\gamma, \varsigma\in \Omega.$

\end{defn}

A morphism $f: (L,\big\{\bullet_{\alpha},\vee_{\alpha,\beta},\{\cdot, \cdot, \cdot\}_{\beta,\gamma},[\cdot, \cdot, \cdot]_{\alpha,\beta,\gamma}\big\}_{\alpha,\beta,\gamma\in \Omega})\rightarrow (L',\big\{\bullet'_{\alpha},\vee'_{\alpha,\beta},\{\cdot, \cdot, \cdot\}'_{\beta,\gamma},[\cdot, \cdot, \cdot]'_{\alpha,\beta,\gamma}\big\}_{\alpha,\beta,\gamma\in \Omega})$ of  NS-Lie-Yamaguti family algebras is a linear
map $f:  L \rightarrow  L' $
satisfying $f(x\bullet_{\alpha} y)=f(x)\bullet'_{\alpha} f(y), f(x\vee_{\alpha,\beta} y)=f(x)\vee'_{\alpha,\beta} f(y)$,  $f\{x, y, z\}_{\beta,\gamma}=\{f(x), f(y), f(z)\}'_{\beta,\gamma}$
as well as  $f[x, y, z]_{\alpha,\beta,\gamma}=[f(x), f(y), f(z)]'_{\alpha,\beta,\gamma}$  for  $x,y,z\in L$ and $\alpha,\beta,\gamma\in \Omega$.

\begin{coro}
We have the following equalities:
 \begin{align}
0=&  \{[x, y]^\ast_{\alpha,\beta},z,a\}_{\alpha\beta,\gamma}^\ast+\{[y, z]_{\beta,\gamma}^\ast,x,a\}_{\beta\gamma,\alpha}^\ast+\{[z, x]_{\gamma,\alpha}^\ast,y,a\}_{\gamma\alpha,\beta}^\ast,\label{4.28}\\
0=&\{x,y,\{z,a,b\}_{\gamma,\varsigma}^\ast\}_{\alpha,\beta}^\ast-\{z,a,\{x,y,b\}_{\alpha,\beta}^\ast\}_{\gamma,\varsigma}^\ast-\{[\![x,y,z]\!]_{\alpha,\beta,\gamma},a,b\}_{\alpha\beta\gamma,\varsigma}^\ast-\label{4.29}\\
&\{z,[\![x,y,a]\!]_{\alpha,\beta,\varsigma},b\}_{\gamma,\alpha\beta\varsigma}^\ast,\nonumber\\
0=&\{b, [\![x,y,z]\!]_{\alpha,\beta,\gamma}, a\}_{\alpha\beta\gamma,\varsigma}-\{\{b, z, y\}_{\gamma,\beta}, x, a\}_{\alpha,\varsigma}+\{\{b, z, x\}_{\gamma,\alpha}, y, a\}_{\beta,\varsigma}+\{\{x, y, b\}^\ast_{\alpha,\beta}, z, a\}_{\gamma,\varsigma}.\label{4.30}
 \end{align}
\end{coro}

\begin{remark}
(i)  Any  NS-Lie-Yamaguti algebra  $(L,\bullet,\vee,\{\cdot, \cdot, \cdot\},[\cdot, \cdot, \cdot])$  can be considered as a constant  NS-pre-Lie family algebra $(L,\big\{\bullet_{\alpha},\vee_{\alpha,\beta},\{\cdot, \cdot, \cdot\}_{\beta,\gamma},[\cdot, \cdot, \cdot]_{\alpha,\beta,\gamma}\big\}_{\alpha,\beta,\gamma\in \Omega})$,
 where $\bullet_{\alpha}=\bullet,\vee_{\alpha,\beta}=\vee,\{\cdot, \cdot, \cdot\}_{\beta,\gamma}=\{\cdot, \cdot, \cdot\}$ and $[\cdot, \cdot, \cdot]_{\alpha,\beta,\gamma}=[\cdot, \cdot, \cdot]$ for all $\alpha,\beta,\gamma\in\Omega$.

(ii) On the one hand, if the   family
of bilinear operations $\{\bullet_{\alpha}\}_{\alpha\in \Omega}$ and  the family
of trilinear operations $\big\{\{\cdot, \cdot, \cdot\}_{\beta,\gamma}\big\}_{\beta,\gamma\in \Omega}$in the above definition are trivial,
 we arrive at the concept of a Lie-Yamaguti family algebra $(L,\big\{\vee_{\alpha,\beta},[\cdot, \cdot, \cdot]_{\alpha,\beta,\gamma}\big\}_{\alpha,\beta,\gamma\in \Omega})$.
On the other hand, if  the families $\{\vee_{\alpha,\beta}\}_{\alpha,\beta\in \Omega}$ and $\{[\cdot, \cdot, \cdot]_{\alpha,\beta,\gamma}\}_{\alpha,\beta,\gamma\in \Omega}$ are trivial,
we obtain the notion of a pre-Lie-Yamaguti family algebra $(L,\big\{\bullet_{\alpha},\{\cdot, \cdot, \cdot\}_{\beta,\gamma}\big\}_{\alpha,\beta,\gamma\in \Omega})$.
\end{remark}

 \begin{prop} \label{prop:4.5}
 Let   $(L,\big\{\bullet_{\alpha},\vee_{\alpha,\beta},\{\cdot, \cdot, \cdot\}_{\beta,\gamma},[\cdot, \cdot, \cdot]_{\alpha,\beta,\gamma}\big\}_{\alpha,\beta,\gamma\in \Omega})$
   be an   NS-Lie-Yamaguti  family algebra.
 Then $(L\otimes \mathbb{K}\Omega,\bullet,\vee,\{\cdot, \cdot, \cdot\},[\cdot, \cdot, \cdot])$  is an  NS-Lie-Yamaguti algebra, where
 \begin{align*}
(x\otimes \alpha)\bullet (y\otimes\beta)=(x\bullet_\alpha y)\otimes\alpha\beta,~\ ~&(x\otimes \alpha)\vee (y\otimes\beta)=(x\vee_{\alpha,\beta} y)\otimes\alpha\beta,\\
\{x\otimes \alpha,   y\otimes\beta, z\otimes \gamma\}=\{x,y,z\}_{\beta,\gamma}\otimes \alpha\beta\gamma,~\ ~&[x\otimes \alpha,   y\otimes\beta, z\otimes \gamma]=[x,y,z]_{\alpha,\beta,\gamma}\otimes \alpha\beta\gamma,
  \end{align*}
for all  $ x,y,z\in L, ~~\alpha,\beta,\gamma\in \Omega.$  Moreover, suppose that $(L',\big\{\bullet'_{\alpha},\vee'_{\alpha,\beta},\{\cdot, \cdot, \cdot\}'_{\beta,\gamma},[\cdot, \cdot, \cdot]'_{\alpha,\beta,\gamma}\big\}_{\alpha,\beta,\gamma\in \Omega})$  is another  NS-Lie-Yamaguti  family algebra.
   If $f:L\rightarrow L'$ is a morphism of  NS-Lie-Yamaguti  family algebras,
  then the map
$\bar{f}:L\otimes \mathbb{K}\Omega\rightarrow L'\otimes \mathbb{K}\Omega, x\otimes \alpha\mapsto f(x)\otimes \alpha$ is a morphism between the induced
NS-Lie-Yamaguti    algebras.
 \end{prop}

 \begin{proof}
Firstly, for all $x,y,z,a,b\in L$ and $\alpha,\beta,\gamma,\varsigma\in\Omega$, we have
 \begin{align*}
&[x, y]^\ast_{\alpha,\beta}\otimes\alpha\beta= (x\bullet_\alpha y )\otimes\alpha\beta-(y\bullet_\beta x\big)\otimes\beta\alpha+(x \vee_{\alpha,\beta} y )\otimes\alpha\beta\\
&=(x\otimes \alpha)\bullet (y\otimes\beta)-(y\otimes\beta)\bullet (x\otimes \alpha)+(x\otimes \alpha)\vee (y\otimes\beta)=[x\otimes \alpha, y\otimes\beta]^\ast,\\
&\{x,y,z\}^\ast_{\alpha,\beta}\otimes\alpha\beta\gamma\\
&=\{z,y,x\}_{\beta,\alpha}\otimes\gamma\beta\alpha-\{z,x,y\}_{\alpha,\beta}\otimes\gamma\alpha\beta+\big(x\bullet_\alpha (y\bullet_\beta z)\big)\otimes\alpha\beta\gamma-\big(y\bullet_\beta (x\bullet_\alpha z)\big)\otimes\beta\alpha\gamma-\big([x,y]^\ast_{\alpha,\beta}\bullet_{\alpha\beta} z\big)\otimes\alpha\beta\gamma\\
&=\{z\otimes\gamma,y\otimes\beta,x\otimes\alpha\}-\{z\otimes\gamma,x\otimes\alpha,y\otimes \beta\}+(x\otimes\alpha)\bullet \big((y\otimes \beta)\bullet (z\otimes\gamma)\big)-\\
&  (y\otimes\beta)\bullet \big((x\otimes \alpha)\bullet (z\otimes\gamma)\big)- [x\otimes\alpha,y\otimes \beta]^\ast  \bullet (z\otimes\gamma)  \\
&=\{x\otimes\alpha,y\otimes \beta,z\otimes\gamma\}^\ast,\\
&[\![x,y,z]\!]_{\alpha,\beta,\gamma}\otimes\alpha\beta\gamma\\
&= \{x,y,z\}^\ast_{\alpha,\beta}\otimes\alpha\beta\gamma+\{x,y,z\}_{\beta,\gamma}\otimes\alpha\beta\gamma-\{y,x,z\}_{\alpha,\gamma}\otimes\beta\alpha\gamma+[x,y,z]_{\alpha,\beta,\gamma}\otimes\alpha\beta\gamma\\
&=\{x\otimes\alpha,y\otimes \beta,z\otimes\gamma\}^\ast+\{x\otimes \alpha,   y\otimes\beta, z\otimes \gamma\}-\{y\otimes\beta, x\otimes \alpha,   z\otimes \gamma\}+[x\otimes \alpha,   y\otimes\beta, z\otimes \gamma]\\
&=[\![x\otimes\alpha,y\otimes \beta,z\otimes\gamma]\!].
 \end{align*}
It follows that
 \begin{align*}
&\{a\otimes\varsigma,[x\otimes\alpha,y\otimes \beta]^\ast , z\otimes\gamma\} -\{(y\otimes \beta)\bullet (a\otimes\varsigma), x\otimes\alpha,z\otimes\gamma\} +\{(x\otimes\alpha)\bullet (a\otimes\varsigma), y\otimes \beta,z\otimes\gamma\}\\
=&\{a,[x,y]^\ast_{\alpha,\beta}, z\}_{\alpha\beta,\gamma}\otimes\varsigma\alpha\beta\gamma-\{y\bullet_\beta a, x,z\}_{\alpha,\gamma}\otimes\beta\varsigma\alpha\gamma+
\{x\bullet_\alpha a, y,z\}_{\beta,\gamma}\otimes\alpha\varsigma\beta\gamma\\
=&0,\\
&\{x\otimes\alpha,y\otimes \beta,(z\otimes\gamma)\bullet (a\otimes\varsigma)\}^\ast-(z\otimes\gamma)\bullet\{x\otimes\alpha,y\otimes \beta,a\otimes\varsigma\}^\ast-[\![x\otimes\alpha,y\otimes \beta,z\otimes\gamma]\!]\bullet (a\otimes\varsigma)\\
=&\{x,y,z\bullet_\gamma a\}^\ast_{\alpha,\beta}\otimes\alpha\beta\gamma\varsigma-z\bullet_\gamma\{x,y,a\}^\ast_{\alpha,\beta}\otimes\gamma\alpha\beta\varsigma-[\![x,y,z]\!]_{\alpha,\beta,\gamma}\bullet_{\alpha\beta\gamma} a\otimes\alpha\beta\gamma\varsigma\\
=&0.
 \end{align*}
Hence Eqs. \eqref{4.1} and \eqref{4.2} hold. Other equations are similarly provable. Therefore, $(L\otimes \mathbb{K}\Omega,\bullet,\vee,\{\cdot, \cdot, \cdot\},[\cdot, \cdot, \cdot])$  is an  NS-Lie-Yamaguti algebra. The last assertion is straightforward
and we complete the proof.
 \end{proof}

 In the following, we give the relations between twisted    Rota-Baxter families and   NS-Lie-Yamaguti  family algebras.

  \begin{prop} \label{prop:4.6}
Let $(V; \rho, \theta)$ be a representation of a   Lie-Yamaguti   algebra $(L, [\cdot, \cdot], \{\cdot, \cdot, \cdot\})$,
 and $\{T_\alpha:V\rightarrow L\}_{\alpha\in\Omega}$ be an $\Gamma$-twisted    Rota-Baxter family.
 Then $(V,\big\{\bullet_{\alpha},\vee_{\alpha,\beta},\{\cdot, \cdot, \cdot\}_{\beta,\gamma},[\cdot, \cdot, \cdot]_{\alpha,\beta,\gamma}\big\}_{\alpha,\beta,\gamma\in \Omega})$ is
  an   NS-Lie-Yamaguti  family algebra, where
 \begin{align*}
u\bullet_{\alpha} v=\rho(T_\alpha u) v,~\ ~ &u\vee_{\alpha,\beta} v= \Gamma_1(T_\alpha u, T_\beta v),\\
\{u, v, w\}_{\beta,\gamma}=\theta(T_\beta v, T_\gamma w)u,~\ ~ &[u, v, w]_{\alpha,\beta,\gamma}= \Gamma_2(T_\alpha u, T_\beta v, T_\gamma w),
  \end{align*}
  for all $u,v,w\in V,~~ \alpha,\beta,\gamma\in \Omega.$
  Furthermore, if $\{T_\alpha:V\rightarrow L\}_{\alpha\in\Omega}$ and $\{T'_\alpha:V'\rightarrow L'\}_{\alpha\in\Omega}$ are two  $\Gamma$-twisted   Rota-Baxter families and $(\zeta,\eta)$ is a morphism between them, then the map $\zeta$ is a morphism between the induced  NS-Lie-Yamaguti  family algebras.
 \end{prop}

\begin{proof}
For any    $u,v,w,s,t\in V$ and $\alpha,\beta,\gamma \in\Omega$. Since
 \begin{align*}
T_{\alpha\beta}[u, v]^\ast_{\alpha,\beta}=&T_{\alpha\beta}(u\bullet_\alpha v-v\bullet_\beta u+u \vee_{\alpha,\beta} v)=T_{\alpha\beta}(\rho(T_\alpha u) v-\rho(T_\beta v) u+ \Gamma_1(T_\alpha u, T_\beta v))\\
=&[T_\alpha u,T_\beta v],\\
 \{u,v,w\}^\ast_{\alpha,\beta}=&\{w,v,u\}_{\beta,\alpha}-\{w,u,v\}_{\alpha,\beta}+u\bullet_\alpha (v\bullet_\beta w)-v\bullet_\beta (u\bullet_\alpha w)-[u,v]^\ast_{\alpha,\beta}\bullet_{\alpha\beta} w\\
 =&\theta(T_\beta v, T_\alpha u)w-\theta(T_\alpha u, T_\beta v)w+\rho(T_\alpha u)  \rho (T_\beta v) w -\rho(T_\beta v)  \rho (T_\alpha u) w -\rho([T_\alpha u,T_\beta v]) w\\
 =&D(T_\alpha u, T_\beta v)w,\\
 T_{\alpha\beta\gamma}[\![u,v,w]\!]_{\alpha,\beta,\gamma}=& T_{\alpha\beta\gamma}\big(\{u,v,w\}^\ast_{\alpha,\beta}+\{u,v,w\}_{\beta,\gamma}-\{v,u,w\}_{\alpha,\gamma}+[u,v,w]_{\alpha,\beta,\gamma}\big)\\
 =&T_{\alpha\beta\gamma}\big(D(T_\alpha u, T_\beta v)w+\theta(T_\beta v, T_\gamma w)u-\theta(T_\alpha u, T_\gamma w)v+\Gamma_2(T_\alpha u, T_\beta v, T_\gamma w)\big)\\
 =&\{T_\alpha u, T_\beta v, T_\gamma w\}
 \end{align*}
and Eqs.\eqref{2.5}-\eqref{2.9},
we have
 \begin{align*}
&\{s, [u, v]^\ast_{\alpha,\beta}, w\}_{\alpha\beta,\gamma}-\{v\bullet_{\beta} s, u, w\}_{\alpha,\gamma}+\{u\bullet_{\alpha} s, v, w\}_{\beta,\gamma}\\
=&\theta([T_\alpha u,T_\beta v],T_\gamma w)s-\theta(T_\alpha u,T_\gamma w)\rho(T_\beta v)s+\theta(T_\beta v,T_\gamma w)\rho(T_\alpha u)s=0,\\
& \{u,v,w\bullet_{\gamma} s\}^\ast_{\alpha,\beta}-w\bullet_{\gamma} \{u,v,s\}^\ast_{\alpha,\beta}-[\![u,v,w]\!]_{\alpha,\beta,\gamma}\bullet_{\alpha\beta\gamma} s\\
=&D(T_\alpha u,T_\beta v)\rho(T_\gamma w)s-\rho(T_\gamma w)D(T_\alpha u,T_\beta v)s-\rho(\{T_\alpha u,T_\beta v,T_\gamma w\})s=0,\\
&\{s, u, [v, w]^\ast_{\beta,\gamma}\}_{\alpha,\beta\gamma}-v\bullet_{\beta} \{s,u,  w\}_{\alpha,\gamma}+w\bullet_{\gamma} \{s,u,  v\}_{\alpha,\beta}\\
=&\theta(T_\alpha u,[T_\beta v,T_\gamma w])s-\rho(T_\beta v)\theta(T_\alpha u,T_\gamma w)s+\rho(T_\gamma w)\theta(T_\alpha u,T_\beta v)s=0,\\
& \{u,v,\{t, w, s\}_{\gamma,\varsigma}\}^\ast_{\alpha,\beta}-\{\{u,v,t\}^\ast_{\alpha,\beta}, w, s\}_{\gamma,\varsigma}-\{t, [\![u,v,w]\!]_{\alpha,\beta,\gamma}, s\}_{\alpha\beta\gamma,\varsigma}-
\{t, w, [\![u,v,s]\!]_{\alpha,\beta,\varsigma}\}_{\gamma,\alpha\beta\varsigma}\\
=&D(T_\alpha u,T_\beta v)\theta(T_\gamma  w,T_\varsigma s)t-\theta(T_\gamma  w,T_\varsigma s)D(T_\alpha u, T_\beta v)t-\theta(\{T_\alpha u,T_\beta v,T_\gamma  w\},T_\varsigma s)t-\\
&\theta(T_\gamma w,\{T_\alpha u,T_\beta v,T_\varsigma s\})t=0,\\
&\{t, u, [\![v,w,s]\!]_{\beta,\gamma,\varsigma}\}_{\alpha,\beta\gamma\varsigma}-\{\{t, u, v\}_{\alpha,\beta}, w, s\}_{\gamma,\varsigma}+\{\{t, u, w\}_{\alpha,\gamma}, v, s\}_{\beta,\varsigma}
- \{v,w,\{t, u, s\}_{\alpha,\varsigma}\}^\ast_{\beta,\gamma}\\
=&\theta(T_\alpha u,\{T_\beta v,T_\gamma  w,T_\varsigma s\})t -\theta(T_\gamma  w,T_\varsigma s)\theta(T_\alpha u,T_\beta v)t+\theta(T_\beta v,T_\varsigma s)\theta(T_\alpha u,T_\gamma  w)t-\\
&D(T_\beta v,T_\gamma  w)\theta(T_\alpha u,T_\varsigma s)t=0.
 \end{align*}
Finally, by $(\Gamma_1,\Gamma_2)$ is a (2,3)-cocycle, we obtain
 \begin{align*}
&[u, v]^\ast_{\alpha,\beta} \vee_{\alpha\beta,\gamma} w+[v, w]^\ast_{\beta,\gamma} \vee_{\beta\gamma,\alpha} u+[w, u]^\ast_{\gamma,\alpha} \vee_{\gamma\alpha,\beta} v-w\bullet_\gamma (u \vee_{\alpha,\beta} v)-u\bullet_\alpha (v \vee_{\beta,\gamma} w)-\\
&v\bullet_\beta (w \vee_{\gamma,\alpha} u)+[u,v,w]_{\alpha,\beta,\gamma}+[v,w,u]_{\beta,\gamma,\alpha}+[w,u,v]_{\gamma,\alpha,\beta}\\
=&\Gamma_1([T_\alpha u,T_\beta v], T_{\gamma} w)+\Gamma_1([T_\beta v,T_\gamma w], T_{\alpha} u)+\Gamma_1([T_\gamma w,T_\alpha u], T_{\beta} v)-
\rho(T_\gamma w)\Gamma_1(T_{\alpha}u, T_{\beta}v)-\\
&\rho(T_\alpha u)\Gamma_1(T_{\beta}v, T_{\gamma}w)-\rho(T_\beta v)\Gamma_1(T_{\gamma}w, T_{\alpha}u)+\Gamma_2(T_\alpha u, T_\beta v, T_\gamma w)+\Gamma_2(T_\beta v, T_\gamma w, T_\alpha u)+\Gamma_2(T_\gamma w, T_\alpha u, T_\beta v)\\
=&0,\\
&\{u \vee_{\alpha,\beta} v,w,s\}_{\gamma,\varsigma}+\{v \vee_{\beta,\gamma} w,u,s\}_{\alpha,\varsigma}+\{w \vee_{\gamma,\alpha} u,v,s\}_{\beta,\varsigma}+[[u, v]^\ast_{\alpha,\beta},w,s]_{\gamma,\varsigma}+ \\
&[[v, w]^\ast_{\beta,\gamma},u,s]_{\alpha,\varsigma}+[[w, u]^\ast_{\gamma,\alpha},v,s]_{\beta,\varsigma}\\
=&\theta(T_\gamma w,T_\varsigma s)\Gamma_1(T_\alpha u, T_\beta v)+\theta(T_\alpha u,T_\varsigma s)\Gamma_1(T_\beta v, T_\gamma w)+\theta(T_\beta v,T_\varsigma s)\Gamma_1(T_\gamma w, T_\alpha u)+ \\
&\Gamma_2([T_\alpha u,T_\beta v],T_\gamma w,T_\varsigma s)+\Gamma_2([T_\beta v,T_\gamma w],T_\alpha u,T_\varsigma s)+\Gamma_2([T_\gamma w,T_\alpha u],T_\beta v,T_\varsigma s)\\
=&0,\\
&-w\bullet_\gamma [u, v, s]_{\alpha,\beta,\varsigma}+s\bullet_\varsigma [u, v, w]_{\alpha,\beta,\gamma}+[u, v, [w, s]^\ast_{\gamma,\varsigma}]_{\alpha,\beta,\gamma}+\{u,v,w\vee_{\gamma,\varsigma} s\}^\ast_{\alpha,\beta}-\\
&[\![u,v,w]\!]_{\alpha,\beta,\gamma}\vee_{\alpha\beta\gamma,\varsigma} s-w\vee_{\gamma,\alpha\beta\varsigma} [\![u,v,s]\!]_{\alpha,\beta,\varsigma}\\
=&-\rho(T_\gamma w)\Gamma_2(T_\alpha u,T_\beta v,T_\varsigma s)+\rho(T_\varsigma s)\Gamma_2(T_\alpha u,T_\beta v,T_\gamma w)+\Gamma_2(T_\alpha u,T_\beta v,[T_\gamma w,T_\varsigma s])+\\
&D(T_\alpha u,T_\beta v)\Gamma_1(T_\gamma w,T_\varsigma s)-\Gamma_1(\{T_\alpha u,T_\beta v,T_\gamma w\},T_\varsigma s)-\Gamma_1(T_\gamma w,\{T_\alpha u,T_\beta v,T_\varsigma s\})\\
=&0
 \end{align*}
 and
 \begin{align*}
& -\{[u, v, w]_{\alpha,\beta,\gamma}, s, t\}_{\varsigma,\tau}+\{[u, v, s]_{\alpha,\beta,\varsigma}, w, t\}_{\gamma,\tau}+\{u,v,[w, s, t]_{\gamma,\varsigma,\tau}\}^\ast_{\alpha,\beta}-\{w,s,[u, v, t]_{\alpha,\beta,\tau}\}^\ast_{\gamma,\varsigma}-\\
&[[\![u,v,w]\!]_{\alpha,\beta,\gamma}, s, t]_{\alpha\beta\gamma,\varsigma,\tau}-[w, [\![u,v,s]\!]_{\alpha,\beta,\varsigma}, t]_{\gamma,\alpha\beta\varsigma,\tau}+[u, v, [\![w,s,t]\!]_{\gamma,\varsigma,\tau}]_{\alpha,\beta,\gamma\varsigma\tau}-[w, s, [\![u,v,t]\!]_{\alpha,\beta,\tau}]_{\gamma,\varsigma,\alpha\beta\tau} \\
=&-\theta(T_\varsigma s,T_\tau t)\Gamma_2(T_\alpha u,T_\beta v,T_\gamma w)+\theta(T_\gamma w,T_\tau t)\Gamma_2(T_\alpha u,T_\beta v,T_\varsigma s)+D(T_\alpha u,T_\beta v)\Gamma_2(T_\gamma w,T_\varsigma s,T_\tau t)-\\
&D(T_\gamma w,T_\varsigma s)\Gamma_2(T_\alpha u,T_\beta v,T_\tau t)-\Gamma_2(\{T_\alpha u,T_\beta v,T_\gamma w\},T_\varsigma s,T_\tau t)-\Gamma_2(T_\gamma w,\{T_\alpha u,T_\beta v,T_\varsigma s\},T_\tau t)+\\
&\Gamma_2(T_\alpha u, T_\beta v, \{T_\gamma w,T_\varsigma s,T_\tau t\})-\Gamma_2(T_\gamma w,T_\varsigma s,\{T_\alpha u,T_\beta v,T_\tau t\})\\
=&0.
 \end{align*}
 Thus, $(V,\big\{\bullet_{\alpha},\vee_{\alpha,\beta},\{\cdot, \cdot, \cdot\}_{\beta,\gamma},[\cdot, \cdot, \cdot]_{\alpha,\beta,\gamma}\big\}_{\alpha,\beta,\gamma\in \Omega})$ is
  an   NS-Lie-Yamaguti  family algebra.
   The last assertion follows from a
routine computation and we omit the details. This completes the proof.
 \end{proof}

 The following corollaries can be obtained   from Proposition \ref{prop:4.6}.

  \begin{coro} \label{coro:4.7}

Let $(L, [\cdot, \cdot], \{\cdot, \cdot, \cdot\})$  be a   Lie-Yamaguti   algebra,  $(V; \rho, \theta)$ be a representation of $L$,  and $(\Gamma_1,\Gamma_2)$  be a (2,3)-cocycle.  If $T:V\rightarrow L$  is a generalized Reynolds operator given in Example \ref{exam:generalized Reynolds operator},
 then $(V, \bullet,\vee,\{\cdot, \cdot, \cdot\},[\cdot, \cdot, \cdot])$ is
  an   NS-Lie-Yamaguti  algebra, where
 \begin{align*}
u\bullet v=\rho(Tu) v,~\ ~ &u\vee v= \Gamma_1(T u, T v),\\
\{u, v, w\}=\theta(T v, T w)u,~\ ~ &[u, v, w]= \Gamma_2(T u, T v, T w),
  \end{align*}
  for all $u,v,w\in V$.
 \end{coro}

  \begin{coro}
Let  $(L, [\cdot, \cdot], \{\cdot, \cdot, \cdot\})$  be a   Lie-Yamaguti   algebra and $\{N_\alpha:L\rightarrow L\}_{\alpha\in\Omega}$  be a    Nijenhuis family  given in Example \ref{exam:Nijenhuis family}.
 Then $(L,\big\{\bullet_{\alpha},\vee_{\alpha,\beta},\{\cdot, \cdot, \cdot\}_{\beta,\gamma},[\cdot, \cdot, \cdot]_{\alpha,\beta,\gamma}\big\}_{\alpha,\beta,\gamma\in \Omega})$ is
  an   NS-Lie-Yamaguti  family algebra, where
 \begin{align*}
x\bullet_{\alpha} y=[N_\alpha x, y],~\ ~ &x\vee_{\alpha,\beta} y=-N_{\alpha\gamma}[x,y],\\
\{x, y, z\}_{\beta,\gamma}=\{z, N_\beta x, N_\gamma y\},~\ ~ &[x, y, z]_{\alpha,\beta,\gamma}= -N_{\alpha\beta\gamma}\big(\{N_\alpha x, y, z\}+\{x,  N_\beta y,z\}+\{x,  y, N_\beta z\}-N_{\alpha\beta\gamma}\{x, y,  z\}\big),
  \end{align*}
  for all $x,y,z\in V,~~ \alpha,\beta,\gamma\in \Omega.$
 \end{coro}

$\bullet $  The first NS-Lie-Yamaguti   algebra is obtained by applying Proposition \ref{prop:4.5} to  NS-Lie-Yamaguti  family algebra
$(V,\big\{\bullet_{\alpha},\vee_{\alpha,\beta},\{\cdot, \cdot, \cdot\}_{\beta,\gamma},[\cdot, \cdot, \cdot]_{\alpha,\beta,\gamma}\big\}_{\alpha,\beta,\gamma\in \Omega})$  defined in Proposition \ref{prop:4.6}. Then
$(V\otimes \mathbb{K}\Omega,\bullet,$ $\vee,\{\cdot, \cdot, \cdot\},[\cdot, \cdot, \cdot])$   is an  NS-Lie-Yamaguti algebra,
where
 \begin{align}
\left\{ \begin{array}{llll}
 &(u\otimes \alpha)\bullet (v\otimes\beta)=(u\bullet_\alpha v)\otimes\alpha\beta=\rho(T_\alpha u)v\otimes\alpha\beta,\\
& (u\otimes \alpha )\vee (v\otimes\beta)=(u \vee_{\alpha,\beta} v)\otimes\alpha\beta=\Gamma_1(T_\alpha u, T_\beta v)\otimes\alpha\beta,\\
& \{u \otimes \alpha, v \otimes \beta, w \otimes \gamma\}=\{u,v,w\}_{\beta,\gamma}\otimes\alpha\beta\gamma=\theta(T_\beta v, T_\gamma w)u \otimes\alpha\beta\gamma,\\
& [u \otimes \alpha, v\otimes\beta, w\otimes \gamma] =[u,v,w]_{\alpha,\beta,\gamma}\otimes\alpha\beta\gamma=\Gamma_2(T_\alpha u, T_\beta v, T_\gamma w)\otimes\alpha\beta\gamma.
\end{array}  \right. \label{4.31}
 \end{align}

$\bullet $   Note that the $\Gamma$-twisted   Rota-Baxter family    $\{T_\alpha:V\rightarrow L\}_{\alpha\in\Omega}$
 induces a generalized Reynolds operator $\overline{T}:V\otimes \mathbb{K}\Omega\rightarrow L\otimes \mathbb{K}\Omega$, $\overline{T}(u\otimes \alpha)=(T_\alpha u)\otimes\alpha$ by Proposition  \ref{prop:2.9}.
 Hence we get the
second NS-Lie-Yamaguti   algebra $V\otimes \mathbb{K}\Omega$ by Corollary \ref{coro:4.7}, which coincides with the one given by \eqref{4.31}.

\section{ $\Omega$-Lie-Yamaguti  algebras and  their cohomology} \label{sec:O-Cohomology}
\def\theequation{\arabic{section}.\arabic{equation}}
\setcounter{equation} {0}

In this section,  we   introduce the notion of $\Omega$-Lie-Yamaguti  algebras,  and  their cohomology

We  assume that $\Omega$ is a  commutative semigroup.

\begin{defn}
 An  $\Omega$-Lie-Yamaguti  algebra consists of a  vector space $L$  equipped with  a collection
of bilinear operations  $\{[\cdot, \cdot]_{\alpha,\beta}:  L\times L\rightarrow L\}_{\alpha,\beta\in \Omega}$ and  a collection
of   trilinear operations   $\big\{\{\cdot, \cdot, \cdot\}_{\alpha,\beta,\gamma}: L \times L\times L\rightarrow L\big\}_{\alpha,\beta,\gamma\in \Omega}$ such that the following conditions hold
\begin{align}
&[x,y]_{\alpha,\beta}=-[y, x]_{\beta, \alpha},\ ~\ ~ \{x,y, z\}_{\alpha,\beta,\gamma}=-\{y, x, z\}_{\beta,\alpha,\gamma},\label{5.1}\\
&[[x, y]_{\alpha,\beta}, z]_{\alpha\beta,\gamma}+[[y, z]_{\beta,\gamma}, x]_{\beta\gamma,\alpha}+[[z, x]_{\gamma,\alpha}, y]_{\gamma\alpha,\beta}+\label{5.2}\\
&\ ~\ ~ \ ~\ ~\ ~\ ~ \ ~\ ~ \ ~\ ~\ ~\ ~  \ ~\ ~ \ ~\ ~\ ~\ ~  \ ~\ ~ \ ~\ ~\ ~\ ~  \ ~\ ~ \ ~\ ~\ ~\ ~  \ ~\ ~ \ ~\ ~ \{x, y, z\}_{\alpha,\beta,\gamma}+\{y, z, x\}_{\beta,\gamma,\alpha}+\{z, x, y\}_{\gamma,\alpha,\beta}=0,\nonumber\\
& \{[x, y]_{\alpha,\beta}, z, a\}_{\alpha\beta,\gamma,\varsigma}+ \{[y, z]_{\beta,\gamma}, x, a\}_{\beta\gamma,\alpha,\varsigma}+ \{[z, x]_{\gamma,\alpha}, y, a\}_{\gamma\alpha,\beta,\varsigma}=0,\label{5.3}\\
&\{a, x, [y, z]_{\beta,\gamma}\}_{\varsigma,\alpha, \beta\gamma}=[\{a, x, y\}_{\varsigma, \alpha, \beta}, z]_{\varsigma\alpha\beta,\gamma}+[y,\{a, x, z\}_{\varsigma,\alpha,\gamma}]_{\beta,\varsigma\alpha\gamma},\label{5.4}\\
& \{a, b, \{x, y, z\}_{\alpha,\beta,\gamma}\}_{\varsigma,\tau,\alpha\beta\gamma}\label{5.5}\\
&=\{\{a, b, x\}_{\varsigma,\tau,\alpha}, y, z\}_{\varsigma\tau\alpha,\beta,\gamma}+ \{x,  \{a, b, y\}_{\varsigma,\tau,\beta}, z\}_{\alpha,\varsigma\tau\beta,\gamma}+ \{x,  y, \{a, b, z\}_{\varsigma,\tau,\gamma}\}_{\alpha,\beta,\varsigma\tau\gamma},\nonumber
\end{align}
for all $ x, y, z, a, b\in L$ and $\alpha,\beta,\gamma,\varsigma,\tau\in \Omega$.
\end{defn}

\begin{exam}
Recall that in \cite{Aguiar2020},  an  $\Omega$-Lie algebra $(L, \{[\cdot, \cdot]_{\alpha,\beta} \}_{\alpha,\beta\in \Omega})$ is defined as a vector space $L$ equipped with an operation $[\cdot, \cdot]_{\alpha,\beta}$ for each pair $(\alpha,\beta)\in \Omega^2$, satisfying the conditions:
\begin{align*}
&[x,y]_{\alpha,\beta}+[y, x]_{\beta, \alpha}=0,\\
&[[x,y]_{\alpha,\beta},z]_{\alpha\beta,\gamma}+[[y,z]_{\beta,\gamma},x]_{\beta\gamma,\alpha}+[[z,x]_{\gamma,\alpha},y]_{\gamma\alpha,\beta}=0,
\end{align*}
for all $ x, y, z\in L$ and $\alpha,\beta,\gamma\in \Omega$.
By defining
\begin{align*}
&\{x,y,z\}_{\alpha,\beta,\gamma}=[[x,y]_{\alpha,\beta},z]_{\alpha\beta,\gamma},
\end{align*}
we obtain the induced $\Omega$-Lie-Yamaguti algebra $(L,\{[\cdot, \cdot]_{\alpha,\beta}, \{\cdot, \cdot, \cdot\}_{\alpha,\beta,\gamma} \}_{\alpha,\beta,\gamma\in \Omega})$.
\end{exam}

   \begin{prop} \label{prop:OLYA}
 Let   $(L,\big\{\bullet_{\alpha},\vee_{\alpha,\beta},\{\cdot, \cdot, \cdot\}_{\beta,\gamma},[\cdot, \cdot, \cdot]_{\alpha,\beta,\gamma}\big\}_{\alpha,\beta,\gamma\in \Omega})$
   be an   NS-Lie-Yamaguti  family algebra.
 Then $(L,\{[\cdot, \cdot]^\ast_{\alpha,\beta}, [\![\cdot, \cdot, \cdot]\!]_{\alpha,\beta,\gamma} \}_{\alpha,\beta,\gamma\in \Omega})$ is an $\Omega$-Lie-Yamaguti  algebra, where
 $[\cdot, \cdot]^\ast_{\alpha,\beta}$ and $[\![\cdot, \cdot, \cdot]\!]_{\alpha,\beta,\gamma}$ are defined by Eqs. \eqref{4.25} and \eqref{4.27}, respectively.
 \end{prop}
    \begin{proof}
For all $ x, y, z, a, b\in L$ and $\alpha,\beta,\gamma,\varsigma,\tau\in \Omega$, it is straightforward to demonstrate that
$[x,y]^\ast_{\alpha,\beta}=-[y, x]^\ast_{\beta, \alpha}$,   $[\![x,y, z]\!]_{\alpha,\beta,\gamma}=-[\![y, x, z]\!]_{\beta,\alpha,\gamma}$,
and by Eq. \eqref{4.21},  we obtain:
\begin{align*}
&[[x, y]^\ast_{\alpha,\beta}, z]^\ast_{\alpha\beta,\gamma}+[[y, z]^\ast_{\beta,\gamma}, x]^\ast_{\beta\gamma,\alpha}+[[z, x]^\ast_{\gamma,\alpha}, y]^\ast_{\gamma\alpha,\beta}+
 [\![x, y, z]\!]_{\alpha,\beta,\gamma}+[\![y, z, x]\!]_{\beta,\gamma,\alpha}+[\![z, x, y]\!]_{\gamma,\alpha,\beta}\\
 =&[x, y]^\ast_{\alpha,\beta}\bullet_{\alpha\beta} z-z\bullet_\gamma (x\bullet_\alpha y-y\bullet_\beta x+x \vee_{\alpha,\beta} y)+[x, y]^\ast_{\alpha,\beta} \vee_{\alpha\beta,\gamma} z+\\
 &[y, z]^\ast_{\beta,\gamma}\bullet_{\beta\gamma} x-x\bullet_\alpha (y\bullet_\beta z-z\bullet_\gamma y+y \vee_{\beta,\gamma} z)+[y, z]^\ast_{\beta,\gamma} \vee_{\beta\gamma,\alpha} x+\\
 &[z, x]^\ast_{\gamma,\alpha}\bullet_{\gamma\alpha} y-y\bullet_\beta (z\bullet_\gamma x-x\bullet_\alpha z+z \vee_{\gamma,\alpha} x)+[z, x]^\ast_{\gamma,\alpha} \vee_{\gamma\alpha,\beta} y+\\
 &\{z,y,x\}_{\beta,\alpha}-\{z,x,y\}_{\alpha,\beta}+x\bullet_\alpha (y\bullet_\beta z)-y\bullet_\beta (x\bullet_\alpha z)-[x,y]^\ast_{\alpha,\beta}\bullet_{\alpha\beta} z+\{x,y,z\}_{\beta,\gamma}-\{y,x,z\}_{\alpha,\gamma}+[x,y,z]_{\alpha,\beta,\gamma}+\\
 &\{x,z,y\}_{\gamma,\beta}-\{x,y,z\}_{\beta,\gamma}+y\bullet_\beta (z\bullet_\gamma x)-z\bullet_\gamma (y\bullet_\beta x)-[y,z]^\ast_{\beta,\gamma}\bullet_{\beta\gamma} x+\{y,z,x\}_{\gamma,\alpha}-\{z,y,x\}_{\beta,\alpha}+[y,z,x]_{\beta,\gamma,\alpha}+\\
 &\{y,x,z\}_{\alpha,\gamma}-\{y,z,x\}_{\gamma,\alpha}+z\bullet_\gamma (x\bullet_\alpha y)-x\bullet_\alpha (z\bullet_\gamma y)-[z,x]^\ast_{\gamma,\alpha}\bullet_{\gamma\alpha} y+\{z,x,y\}_{\alpha,\beta}-\{x,z,y\}_{\gamma,\beta}+[z,x,y]_{\gamma,\alpha,\beta}\\
  =&-z\bullet_\gamma (x \vee_{\alpha,\beta} y)+[x, y]^\ast_{\alpha,\beta} \vee_{\alpha\beta,\gamma} z-x\bullet_\alpha (y \vee_{\beta,\gamma} z)+[y, z]^\ast_{\beta,\gamma} \vee_{\beta\gamma,\alpha} x-y\bullet_\beta (z \vee_{\gamma,\alpha} x)+\\
  &[z, x]^\ast_{\gamma,\alpha} \vee_{\gamma\alpha,\beta} y+[x,y,z]_{\alpha,\beta,\gamma}+[y,z,x]_{\beta,\gamma,\alpha}+[z,x,y]_{\gamma,\alpha,\beta}\\
  =&0.
  \end{align*}
By Eqs. \eqref{4.22} and \eqref{4.28}, we have:
\begin{align*}
 & [\![[x, y]^\ast_{\alpha,\beta}, z, a]\!]_{\alpha\beta,\gamma,\varsigma}+ [\![[y, z]^\ast_{\beta,\gamma}, x, a]\!]_{\beta\gamma,\alpha,\varsigma}+ [\![[z, x]^\ast_{\gamma,\alpha}, y, a]\!]_{\gamma\alpha,\beta,\varsigma}\\
 =& \{[x, y]^\ast_{\alpha,\beta},z,a\}^\ast_{\alpha\beta,\gamma}+\{[x, y]^\ast_{\alpha,\beta},z,a\}_{\gamma,\varsigma}-\{z,[x, y]^\ast_{\alpha,\beta},a\}_{\alpha\beta,\varsigma}+[[x, y]^\ast_{\alpha,\beta},z,a]_{\gamma,\varsigma}+\\
&\{[y, z]^\ast_{\beta,\gamma},x,a\}^\ast_{\beta\gamma,\alpha}+\{[y, z]^\ast_{\beta,\gamma},x,a\}_{\alpha,\varsigma}-\{x,[y, z]^\ast_{\beta,\gamma},a\}_{\beta\gamma,\varsigma}+
[[y, z]^\ast_{\beta,\gamma},x,a]_{\alpha,\varsigma}+\\
&\{[z, x]^\ast_{\gamma,\alpha},y,a\}^\ast_{\gamma\alpha,\beta}+\{[z, x]^\ast_{\gamma,\alpha},y,a\}_{\beta,\varsigma}-\{y,[z, x]^\ast_{\gamma,\alpha},a\}_{\gamma\alpha,\varsigma}+
[[z, x]^\ast_{\gamma,\alpha},y,a]_{\beta,\varsigma}\\
 =& \big\{x\bullet_\alpha y-y\bullet_\beta x+x \vee_{\alpha,\beta} y,z,a\big\}_{\gamma,\varsigma}-\{y\bullet_\beta z, x,a\}_{\alpha,\varsigma}+\{x\bullet_\alpha z, y,a\}_{\beta,\varsigma}+[[x, y]^\ast_{\alpha,\beta},z,a]_{\gamma,\varsigma}+\\
&\big\{y\bullet_\beta z-z\bullet_\gamma y+y \vee_{\beta,\gamma} z,x,a\big\}_{\alpha,\varsigma}-\{z\bullet_\gamma x, y,a\}_{\beta,\varsigma}+\{y\bullet_\beta x, z,a\}_{\gamma,\varsigma}+[[y, z]^\ast_{\beta,\gamma},x,a]_{\alpha,\varsigma}+\\
&\big\{z\bullet_\gamma x-x\bullet_\alpha z+z \vee_{\gamma,\alpha} x,y,a\big\}_{\beta,\varsigma}-\{x\bullet_\alpha y, z,a\}_{\gamma,\varsigma}+\{z\bullet_\gamma y, x,a\}_{\alpha,\varsigma}+[[z, x]^\ast_{\gamma,\alpha},y,a]_{\beta,\varsigma}\\
=& \{x \vee_{\alpha,\beta} y,z,a\}_{\gamma,\varsigma}+[[x, y]^\ast_{\alpha,\beta},z,a]_{\gamma,\varsigma}+\{y \vee_{\beta,\gamma} z,x,a\}_{\alpha,\varsigma}+[[y, z]^\ast_{\beta,\gamma},x,a]_{\alpha,\varsigma}+\\
&\{z \vee_{\gamma,\alpha} x,y,a\}_{\beta,\varsigma}+[[z, x]^\ast_{\gamma,\alpha},y,a]_{\beta,\varsigma}\\
=&0.
\end{align*}
Furthermore, by Eqs. \eqref{4.17}, \eqref{4.18} and \eqref{4.23}, we have
\begin{align*}
&[\![a, x, [y, z]^\ast_{\beta,\gamma}]\!]_{\varsigma,\alpha, \beta\gamma}-[[\![a, x, y]\!]_{\varsigma, \alpha, \beta}, z]^\ast_{\varsigma\alpha\beta,\gamma}-
[y,[\![a, x, z]\!]_{\varsigma,\alpha,\gamma}]^\ast_{\beta,\varsigma\alpha\gamma}\\
=&\{a,x,y\bullet_\beta z\}^\ast_{\varsigma,\alpha}-\{a,x,z\bullet_\gamma y\}^\ast_{\varsigma,\alpha}+\{a,x,y \vee_{\beta,\gamma} z\}^\ast_{\varsigma,\alpha}+\{a,x,[y, z]^\ast_{\beta,\gamma}\}_{\alpha,\beta\gamma}-\{x,a,[y, z]^\ast_{\beta,\gamma}\}_{\varsigma,\beta\gamma}+\\
&[a,x,[y, z]^\ast_{\beta,\gamma}]_{\varsigma,\alpha,\beta\gamma}-[\![a, x, y]\!]_{\varsigma, \alpha, \beta}\bullet_{\varsigma\alpha\beta} z+z\bullet_\gamma \{a,x,y\}^\ast_{\varsigma,\alpha}+z\bullet_\gamma\{a,x,y\}_{\alpha,\beta}-z\bullet_\gamma\{x,a,y\}_{\varsigma,\beta}+\\
&z\bullet_\gamma[a,x,y]_{\varsigma,\alpha,\beta}-[\![a, x, y]\!]_{\varsigma, \alpha, \beta} \vee_{\varsigma\alpha\beta,\gamma} z-y\bullet_{\beta} \{a,x,z\}^\ast_{\varsigma,\alpha}-y\bullet_{\beta}\{a,x,z\}_{\alpha,\gamma}+y\bullet_{\beta}\{x,a,z\}_{\varsigma,\gamma}-\\
&y\bullet_{\beta}[a,x,z]_{\varsigma,\alpha,\gamma}+[\![a, x, z]\!]_{\varsigma,\alpha,\gamma}\bullet_{\varsigma\alpha\gamma} y-y \vee_{\beta,\varsigma\alpha\gamma}[\![a, x, z]\!]_{\varsigma,\alpha,\gamma}\\
=&\{a,x,y \vee_{\beta,\gamma} z\}^\ast_{\varsigma,\alpha}+[a,x,[y, z]^\ast_{\beta,\gamma}]_{\varsigma,\alpha,\beta\gamma} +z\bullet_\gamma[a,x,y]_{\varsigma,\alpha,\beta}-[\![a, x, y]\!]_{\varsigma, \alpha, \beta} \vee_{\varsigma\alpha\beta,\gamma} z-\\
&y\bullet_{\beta}[a,x,z]_{\varsigma,\alpha,\gamma}-y \vee_{\beta,\varsigma\alpha\gamma}[\![a, x, z]\!]_{\varsigma,\alpha,\gamma}\\
=&0.
\end{align*}
  Finally, by Eqs. \eqref{4.19}, \eqref{4.20}, \eqref{4.24} and \eqref{4.29}, we have
  \begin{align*}
&[\![a, b, [\![x, y, z]\!]_{\alpha,\beta,\gamma}]\!]_{\varsigma,\tau,\alpha\beta\gamma}-[\![[\![a, b, x]\!]_{\varsigma,\tau,\alpha}, y, z]\!]_{\varsigma\tau\alpha,\beta,\gamma}-[\![x,  [\![a, b, y]\!]_{\varsigma,\tau,\beta}, z]\!]_{\alpha,\varsigma\tau\beta,\gamma}- [\![x,  y, [\![a, b, z]\!]_{\varsigma,\tau,\gamma}]\!]_{\alpha,\beta,\varsigma\tau\gamma}\\
=&\{a,b,\{x,y,z\}^\ast_{\alpha,\beta}\}^\ast_{\varsigma,\tau}+\{a,b,\{x,y,z\}_{\beta,\gamma}\}^\ast_{\varsigma,\tau}-\{a,b,\{y,x,z\}_{\alpha,\gamma}\}^\ast_{\varsigma,\tau}+\{a,b,[x,y,z]_{\alpha,\beta,\gamma}\}^\ast_{\varsigma,\tau}+\\
&\{a,b,[\![x, y, z]\!]_{\alpha,\beta,\gamma}\}_{\tau,\alpha\beta\gamma}-\{b,a,[\![x, y, z]\!]_{\alpha,\beta,\gamma}\}_{\varsigma,\alpha\beta\gamma}+[a,b,[\![x, y, z]\!]_{\alpha,\beta,\gamma}]_{\varsigma,\tau,\alpha\beta\gamma}-\\
&\{[\![a, b, x]\!]_{\varsigma,\tau,\alpha},y,z\}^\ast_{\varsigma\tau\alpha,\beta}-\{[\![a, b, x]\!]_{\varsigma,\tau,\alpha},y,z\}_{\beta,\gamma}+\{y,[\![a, b, x]\!]_{\varsigma,\tau,\alpha},z\}_{\varsigma\tau\alpha,\gamma}-[[\![a, b, x]\!]_{\varsigma,\tau,\alpha},y,z]_{\varsigma\tau\alpha,\beta,\gamma}-\\
&\{x,[\![a, b, y]\!]_{\varsigma,\tau,\beta},z\}^\ast_{\alpha,\varsigma\tau\beta}-\{x,[\![a, b, y]\!]_{\varsigma,\tau,\beta},z\}_{\varsigma\tau\beta,\gamma}+\{[\![a, b, y]\!]_{\varsigma,\tau,\beta},x,z\}_{\alpha,\gamma}-[x,[\![a, b, y]\!]_{\varsigma,\tau,\beta},z]_{\alpha,\varsigma\tau\beta,\gamma}-\\
&\{x,y,\{a,b,z\}^\ast_{\varsigma,\tau}\}^\ast_{\alpha,\beta}-\{x,y,\{a,b,z\}_{\tau,\gamma}\}^\ast_{\alpha,\beta}+\{x,y,\{b,a,z\}_{\varsigma,\gamma}\}^\ast_{\alpha,\beta}-\{x,y,[a,b,z]_{\varsigma,\tau,\gamma}\}^\ast_{\alpha,\beta}-\\
&\{x,y,[\![a, b, z]\!]_{\varsigma,\tau,\gamma}\}_{\beta,\varsigma\tau\gamma}+\{y,x,[\![a, b, z]\!]_{\varsigma,\tau,\gamma}\}_{\alpha,\varsigma\tau\gamma}-[x,y,[\![a, b, z]\!]_{\varsigma,\tau,\gamma}]_{\alpha,\beta,\varsigma\tau\gamma}\\
=&\{a,b,[x,y,z]_{\alpha,\beta,\gamma}\}^\ast_{\varsigma,\tau}+[a,b,[\![x, y, z]\!]_{\alpha,\beta,\gamma}]_{\varsigma,\tau,\alpha\beta\gamma}-\{[a,b,x]_{\varsigma,\tau,\alpha},y,z\}_{\beta,\gamma}-[[\![a, b, x]\!]_{\varsigma,\tau,\alpha},y,z]_{\varsigma\tau\alpha,\beta,\gamma}+\\
&\{\{a,b,y\}_{\tau,\beta},x,z\}_{\alpha,\gamma}-[x,[\![a, b, y]\!]_{\varsigma,\tau,\beta},z]_{\alpha,\varsigma\tau\beta,\gamma}-\{x,y,[a,b,z]_{\varsigma,\tau,\gamma}\}^\ast_{\alpha,\beta}-[x,y,[\![a, b, z]\!]_{\varsigma,\tau,\gamma}]_{\alpha,\beta,\varsigma\tau\gamma}\\
=&0.
\end{align*}
Therefore, we deduce that $(L,\{[\cdot, \cdot]^\ast_{\alpha,\beta}, [\![\cdot, \cdot, \cdot]\!]_{\alpha,\beta,\gamma} \}_{\alpha,\beta,\gamma\in \Omega})$ is an $\Omega$-Lie-Yamaguti  algebra.
  \end{proof}

  Combining Propositions \ref{prop:4.6} and \ref{prop:OLYA}, we get the following result.

  \begin{coro} \label{coro:O-LY}
Let $(V; \rho, \theta)$ be a representation of a   Lie-Yamaguti   algebra $(L, [\cdot, \cdot], \{\cdot, \cdot, \cdot\})$,
 and $\{T_\alpha:V\rightarrow L\}_{\alpha\in\Omega}$ be an $\Gamma$-twisted    Rota-Baxter family.
 Then $(V,\{[\cdot, \cdot]_{\alpha,\beta}, \{\cdot, \cdot, \cdot\}_{\alpha,\beta,\gamma} \}_{\alpha,\beta,\gamma\in \Omega})$ is
  an $\Omega$-Lie-Yamaguti    algebra, where
 \begin{align*}
[u, v]_{\alpha,\beta}=& \rho(T_\alpha u) v-\rho(T_\beta v) u+ \Gamma_1(T_\alpha u, T_\beta v),\\
\{u,v,w\}_{\alpha,\beta,\gamma}=& D(T_\alpha u, T_\beta v)w+\theta(T_\beta v, T_\gamma w)u-\theta(T_\alpha u, T_\gamma w)v+\Gamma_2(T_\alpha u, T_\beta v, T_\gamma w),
  \end{align*}
  for all $u,v,w\in V$ and $\alpha,\beta,\gamma\in \Omega.$
  \end{coro}

  Moreover, by Example  \ref{exam:Reynolds family}  and Corollary \ref{coro:O-LY}, we have the following result.
    \begin{coro} \label{coro:R-O-LY}
Let $(L, [\cdot, \cdot], \{\cdot, \cdot, \cdot\})$ be a    Lie-Yamaguti   algebra,
 and $\{T_\alpha:L\rightarrow L\}_{\alpha\in\Omega}$ be a Reynolds family.
 Then $(L,\{[\cdot, \cdot]_{\alpha,\beta}, \{\cdot, \cdot, \cdot\}_{\alpha,\beta,\gamma} \}_{\alpha,\beta,\gamma\in \Omega})$ is
  an $\Omega$-Lie-Yamaguti    algebra, where
 \begin{align*}
[x, y]_{\alpha,\beta}=& [T_\alpha x,y]+[x,T_\beta y]-[T_\alpha x,T_\beta y],\\
\{x,y,z\}_{\alpha,\beta,\gamma}=& \{T_\alpha x, T_\beta y, z\}+\{x,T_\beta y, T_\gamma z\}+\{T_\alpha x, y, T_\gamma z\}-\{T_\alpha x, T_\beta y, T_\gamma z\},
  \end{align*}
  for all $x,y,z\in L$ and $\alpha,\beta,\gamma\in \Omega.$
  \end{coro}

\begin{defn}
 Let  $(L,\{[\cdot, \cdot]_{\alpha,\beta}, \{\cdot, \cdot, \cdot\}_{\alpha,\beta,\gamma} \}_{\alpha,\beta,\gamma\in \Omega})$  be  an $\Omega$-Lie-Yamaguti  algebra.  A representation over it
consists of a vector space $V$ together with a collection
  \begin{align*}
  \left\{\begin{array}{cccc}
 \rho_{\alpha,\varsigma}:&L\times V\rightarrow V,~~~ (x,u)\mapsto \rho_{\alpha,\varsigma}(x,u),\\
\theta_{\alpha,\beta,\varsigma}:&L \times L \times V\rightarrow V, ~~~(x,y,u)\mapsto \theta_{\alpha,\beta,\varsigma}(x,y,u)
 \end{array}\right\}_{\alpha,\beta,\varsigma\in\Omega}
 \end{align*}
of bilinear maps satisfying, for $x,y,z,a\in L, u\in V$ and $\alpha,\beta,\gamma,\varsigma,\tau\in \Omega$,
\begin{align}
0=&\theta_{\alpha\beta,\gamma,\varsigma}([x,y]_{\alpha,\beta},z,u)-\theta_{\alpha,\gamma,\beta\varsigma}(x,z,\rho_{\beta,\varsigma}(y,u))+\theta_{\beta,\gamma,\alpha\varsigma}(y,z,\rho_{\alpha,\varsigma}(x,u)),\label{5.6}\\
0=&D_{\alpha,\beta,\gamma\varsigma}(x,y,\rho_{\gamma,\varsigma}(z,u))-\rho_{\gamma,\alpha\beta\varsigma}(z,D_{\alpha,\beta,\varsigma}(x,y,u))-\rho_{\alpha\beta\gamma,\varsigma}(\{x,y,z\}_{\alpha,\beta,\gamma},u),\label{5.7}\\
0=&\theta_{\alpha,\beta\gamma,\varsigma}(x,[y,z]_{\beta,\gamma},u)-\rho_{\beta,\alpha\gamma\varsigma}(y,\theta_{\alpha,\gamma,\varsigma}(x,z,u))+\rho_{\gamma,\alpha\beta\varsigma}(z,\theta_{\alpha,\beta,\varsigma}(x,y,u)),\label{5.8}\\
0=&D_{\tau,\alpha,\beta\gamma\varsigma}(a,x,\theta_{\beta,\gamma,\varsigma}(y,z,u))-\theta_{\beta,\gamma,\tau\alpha\varsigma}(y,z,D_{\tau,\alpha,\varsigma}(a,x,u))-\label{5.9}\\
&\theta_{\tau\alpha\beta,\gamma,\varsigma}(\{a,x,y\}_{\varsigma,\alpha,\beta},z,u)-\theta_{\beta,\tau\alpha\gamma,\varsigma}(y,\{a,x,z\}_{\tau,\alpha,\gamma},u),\nonumber\\
0=&\theta_{\tau,\alpha\beta\gamma,\varsigma}(a,\{x,y,z\}_{\alpha,\beta,\gamma},u)-\theta_{\beta,\gamma,\tau\alpha\varsigma}(y,z,\theta_{\tau,\alpha,\varsigma}(a,x,u))+\label{5.10}\\
&\theta_{\alpha,\gamma,\tau\beta\varsigma}(x,z,\theta_{\tau,\beta,\varsigma}(a,y,u))-D_{\alpha,\beta,\tau\gamma\varsigma}(x,y,\theta_{\tau,\gamma,\varsigma}(a,z,u)).\nonumber
\end{align}
where $\{D_{\alpha,\beta,\varsigma}: L \times L \times V\rightarrow V, ~~~(x,y,u)\mapsto D_{\alpha,\beta,\varsigma}(x,y,u)\}_{\alpha,\beta,\varsigma\in\Omega}$ is given by
\begin{align}
D_{\alpha,\beta,\varsigma}(x,y,u)=&\theta_{\beta,\alpha,\varsigma}(y,x,u)-\theta_{\alpha,\beta,\varsigma}(x,y,u)-\rho_{\alpha\beta,\varsigma}([x,y]_{\alpha,\beta},u)+\label{5.11}\\
&\rho_{\alpha,\beta\varsigma}(x,\rho_{\beta,\varsigma}(y,u))-\rho_{\beta,\alpha\varsigma}(y,\rho_{\alpha,\varsigma}(x,u)).\nonumber
\end{align}
 \end{defn}


Now, we examine the cohomology of an $\Omega$-Lie-Yamaguti  algebra. Let $(V; \{\rho_{\alpha,\varsigma},\theta_{\alpha,\beta,\varsigma}\}_{\alpha,\beta,\varsigma\in\Omega})$ be a
representation of an $\Omega$-Lie-Yamaguti algebra  $(L,$ $\{[\cdot, \cdot]_{\alpha,\beta}, \{\cdot, \cdot, \cdot\}_{\alpha,\beta,\gamma} \}_{\alpha,\beta,\gamma\in \Omega})$.
For $n\geq 1$, we define the space $C^{2n}(L,V)$ as the set of all multilinear maps of the form
$$f_{\alpha_1, \ldots,\alpha_{2n}}: L^{\otimes2n}\rightarrow V, x_1\otimes\cdots\otimes x_{2n}\mapsto f_{\alpha_1,\alpha_2,\ldots,\alpha_{2n}}(x_1, \ldots,x_{2n}) \ ~ \text{for} \ ~  \alpha_1, \ldots,\alpha_{2n} \in\Omega$$
satisfying the following   condition:
$$f_{\alpha_1,\ldots,\alpha_{2i-1},\alpha_{2i},\ldots,\alpha_{2n}}(x_1,\ldots,x_{2i-1},x_{2i},\ldots,x_{2n})+
f_{\alpha_1,\ldots,\alpha_{2i},\alpha_{2i-1},\ldots,\alpha_{2n}}(x_1,\ldots,x_{2i},x_{2i-1},\ldots,x_{2n})=0,$$
for $i=1,\ldots,n.$ Also, we define the space $C^{2n+1}(L,V)$ as the set of all multilinear maps of the form
$$f_{\alpha_1, \ldots,\alpha_{2n+1}}: L^{\otimes2n+1}\rightarrow V, x_1\otimes\cdots\otimes x_{2n+1}\mapsto f_{\alpha_1,\alpha_2,\ldots,\alpha_{2n+1}}(x_1, \ldots,x_{2n+1}) $$
satisfying the following   condition:
$$f_{\alpha_1,\ldots,\alpha_{2i-1},\alpha_{2i},\ldots,\alpha_{2n+1}}(x_1,\ldots,x_{2i-1},x_{2i},\ldots,x_{2n+1})+
f_{\alpha_1,\ldots,\alpha_{2i},\alpha_{2i-1},\ldots,\alpha_{2n+1}}(x_1,\ldots,x_{2i},x_{2i-1},\ldots,x_{2n+1})=0.$$
We define $C^{(2n,2n+1) }_{\mathrm{\Omega}-\mathrm{LY}}(L,V)$,   and $C^{1}_{\mathrm{\Omega}-\mathrm{LY}}(L,V)$ 
 as follows:
\begin{align*}
\left\{ \begin{array}{lll}
C^{(2n,2n+1) }_{\mathrm{\Omega}-\mathrm{LY}}(L,V):=C^{2n}(L,V)\times C^{2n+1}(L,V)\mbox{ \mbox{}  $n\geq 1,$  }\\
C^{1 }_{\mathrm{\Omega}-\mathrm{LY}}(L,V):=\big\{f=\{f_{\alpha}\}_{\alpha\in \Omega}~|~f_{\alpha}: L\rightarrow V~~\text{is linear}\big\}. \mbox{ \mbox{}}
\end{array}  \right.
\end{align*}
There exists   maps are defined as follows:\\
For $n\geq 1$, the map $\delta_{\Omega}=(\delta_{\Omega,\mathrm{I}},\delta_{\Omega,\mathrm{II}}):C^{(2n,2n+1) }_{\mathrm{\Omega}-\mathrm{LY}}(L,V)\rightarrow C^{(2n+2,2n+3) }_{\mathrm{\Omega}-\mathrm{LY}}(L,V)$ is defined by:
\begin{align*}
&(\delta_{\Omega,\mathrm{I}}(f,g))_{\alpha_1,\ldots,\alpha_{2n+2}}(x_1,\ldots,x_{2n+2})\\
=&(-1)^n \Big(\rho_{\alpha_{2n+1},\alpha_1\cdots\alpha_{2n}\alpha_{2n+2}}\big(x_{2n+1}, g_{\alpha_1,\ldots,\alpha_{2n},\alpha_{2n+2}}(x_1,\ldots,x_{2n},x_{2n+2})\big)-\rho_{\alpha_{2n+2},\alpha_1\cdots \alpha_{2n+1}}\big(x_{2n+2}, \\
&g_{\alpha_1,\ldots, \alpha_{2n+1}}(x_1,\ldots, x_{2n+1})\big)-g_{\alpha_1,\ldots, \alpha_{2n},\alpha_{2n+1}\alpha_{2n+2}}\big(x_1,\ldots, x_{2n},[x_{2n+1},x_{2n+2}]_{\alpha_{2n+1},\alpha_{2n+2}}\big)\Big)+\\
&\sum_{k=1}^n(-1)^{k+1}D_{\alpha_{2k-1},\alpha_{2k},\alpha_1\cdots\alpha_{2k-2}\alpha_{2k+1}\cdots\alpha_{2n+2}}\big(x_{2k-1},x_{2k},f_{\alpha_1,\ldots, \alpha_{2k-2},\alpha_{2k+1},\ldots, \alpha_{2n+2}}(x_1,\ldots,\widehat{x}_{2k-1},\widehat{x}_{2k}, \ldots, x_{2n+2})\big)+\\
&\sum_{k=1}^n\sum_{j=2k+1}^{2n+2}(-1)^{k}f_{\alpha_1,\ldots, \alpha_{2k-2},\alpha_{2k+1},\ldots,\alpha_{2k-1}\alpha_{2k}\alpha_{j},\ldots, \alpha_{2n+2}}
\big(x_1,\ldots,\widehat{x}_{2k-1},\widehat{x}_{2k}, \ldots, \{{x}_{2k-1}, {x}_{2k},x_j\}_{\alpha_{2k-1}, \alpha_{2k},\alpha_j}, \ldots,x_{2n+2}\big);\\
&(\delta_{\Omega,\mathrm{II}}(f,g))_{\alpha_1,\ldots,\alpha_{2n+3}}(x_1,\ldots,x_{2n+3})\\
=&(-1)^n \Big(\theta_{\alpha_{2n+2},\alpha_{2n+3},\alpha_1\cdots\alpha_{2n+1}}\big(x_{2n+2},x_{2n+3}, g_{\alpha_1,\ldots, \alpha_{2n+1}}(x_1,\ldots, x_{2n+1})\big)-\\
&\theta_{\alpha_{2n+1},\alpha_{2n+3},\alpha_1\cdots \alpha_{2n}\alpha_{2n+2}}\big(x_{2n+1}, x_{2n+3},g_{\alpha_1,\ldots, \alpha_{2n},\alpha_{2n+2}}(x_1,\ldots, x_{2n}, x_{2n+2})\big)\Big)+\\
&\sum_{k=1}^{n+1}(-1)^{k+1}D_{\alpha_{2k-1},\alpha_{2k},\alpha_1\cdots\alpha_{2k-2}\alpha_{2k+1}\cdots\alpha_{2n+3}}\big(x_{2k-1},x_{2k},g_{\alpha_1,\ldots, \alpha_{2k-2},\alpha_{2k+1},\ldots, \alpha_{2n+3}}(x_1,\ldots,\widehat{x}_{2k-1},\widehat{x}_{2k}, \ldots, x_{2n+3})\big)+\\
&\sum_{k=1}^{n+1}\sum_{j=2k+1}^{2n+3}(-1)^{k}g_{\alpha_1,\ldots, \alpha_{2k-2},\alpha_{2k+1},\ldots,\alpha_{2k-1}\alpha_{2k}\alpha_{j},\ldots, \alpha_{2n+3}}
(x_1,\ldots,\widehat{x}_{2k-1},\widehat{x}_{2k}, \ldots, \{{x}_{2k-1}, {x}_{2k},x_j\}_{\alpha_{2k-1}, \alpha_{2k},\alpha_j}, \ldots,x_{2n+3})\big).
\end{align*}
The map $\delta_{\Omega}=(\delta_{\Omega,\mathrm{I}},\delta_{\Omega,\mathrm{II}}):C^{1 }_{\mathrm{\Omega}-\mathrm{LY}}(L,V)\rightarrow C^{(2,3) }_{\mathrm{\Omega}-\mathrm{LY}}(L,V)$ is defined by:
\begin{align*}
(\delta_{\Omega,\mathrm{I}}f)_{\alpha,\beta}(x,y)=&\rho_{\alpha,\beta}(x,f_{\beta}(y))-\rho_{\beta,\alpha}(y,f_\alpha(x))-f_{\alpha\beta}([x,y]_{\alpha,\beta}),\\
(\delta_{\Omega,\mathrm{II}}f)_{\alpha,\beta,\gamma}(x,y,z)=&D_{\alpha,\beta,\gamma}(x,y,f_\gamma(z))+\theta_{\beta,\gamma,\alpha}(y,z,f_\alpha(x))-
\theta_{\alpha,\gamma,\beta}(x,z,f_{\beta}(y))-f_{\alpha\beta\gamma}\{x,y,z\}_{\alpha,\beta,\gamma}.
\end{align*}
The map $\delta^*_{\Omega}=(\delta^*_{\Omega,\mathrm{I}},\delta^*_{\Omega,\mathrm{II}}):C^{(2,3) }_{\mathrm{\Omega}-\mathrm{LY}}(L,V)\rightarrow C^{(3,4) }_{\mathrm{\Omega}-\mathrm{LY}}(L,V):=C^{3}(L,V)\times C^{4}(L,V)$ is defined by:
\begin{align*}
(\delta^*_{\Omega,\mathrm{I}}(f,g))_{\alpha,\beta,\gamma}(x,y,z)=&-\rho_{\alpha,\beta\gamma}(x,f_{\beta,\gamma}(y,z))-\rho_{\beta,\gamma\alpha}(y,f_{\gamma,\alpha}(z,x))-\rho_{\gamma,\alpha\beta}(z,f_{\alpha,\beta}(x,y))+\\
&f_{\alpha\beta,\gamma}([x,y]_{\alpha,\beta},z)+f_{\beta\gamma,\alpha}([y,z]_{\beta,\gamma},x)+f_{\gamma\alpha,\beta}([z,x]_{\gamma,\alpha},y)+\\
&g_{\alpha,\beta,\gamma}(x,y,z)+g_{\beta,\gamma,\alpha}(y,z,x)+g_{\gamma,\alpha,\beta}(z,x,y),\\
(\delta^*_{\Omega,\mathrm{II}}(f,g))_{\alpha,\beta,\gamma,\varsigma}(x,y,z,a)=&\theta_{\alpha,\varsigma,\beta\gamma}(x,a,f_{\beta,\gamma}(y,z))+\theta_{\beta,\varsigma,\gamma\alpha}(y,a,f_{\gamma,\alpha}(z,x))+
\theta_{\gamma,\varsigma,\alpha\beta}(z,a,f_{\alpha,\beta}(x,y))+\\
&g_{\alpha\beta,\gamma,\varsigma}([x,y]_{\alpha,\beta},z,a)+g_{\beta\gamma,\alpha,\varsigma}([y,z]_{\beta,\gamma},x,a)+g_{\gamma\alpha,\beta,\varsigma}([z,x]_{\gamma,\alpha},y,a).
\end{align*}

  One can verify that $\delta_\Omega \circ \delta_\Omega=0$ by using the same argument as in the usual Lie-Yamaguti  algebra.
  That is, a cochain complex
defined as follows:
$$\aligned
\xymatrix{
 C^{1}_{\Omega-\mathrm{LY}}(L, V)\ar[r]^-{\delta_{\Omega}} &  C^{(2,3)}_{\Omega-\mathrm{LY}}(L,V)\ar[r]^-{\delta_\Omega} \ar[d]^-{\delta_\Omega^*}& C^{(4,5)}_{\Omega-\mathrm{LY}}(L,V)\ar[r]^-{\delta_\Omega} & C^{(6,7)}_{\Omega-\mathrm{LY}}(L,V)\ar[r]^-{\delta_\Omega}&\cdots\\
 ~~& C^{(3,4)}_{\Omega-\mathrm{LY}}(L,V) &  ~& ~&}
 \endaligned$$
  The corresponding cohomology groups are called the cohomology of the $\Omega$-Lie-Yamaguti  algebra
$L$ with coefficients in $V$.

\section{Cohomology and  deformations    of  twisted  Rota-Baxter families} \label{sec:Cohomology}
\def\theequation{\arabic{section}.\arabic{equation}}
\setcounter{equation} {0}

Let $(L, [\cdot, \cdot], \{\cdot, \cdot, \cdot\})$ be a   Lie-Yamaguti   algebra,
$(V; \rho, \theta)$ be a representation of it, and $\Gamma=(\Gamma_1,\Gamma_2)$ be a (2,3)-cocycle.
 Let $\{T_\alpha:V\rightarrow L\}_{\alpha\in\Omega}$ be an $\Gamma$-twisted    Rota-Baxter family.
 Then we have seen in Corollary  \ref{coro:O-LY}
that $(V,\{[\cdot, \cdot]_{\alpha,\beta}, \{\cdot, \cdot, \cdot\}_{\alpha,\beta,\gamma} \}_{\alpha,\beta,\gamma\in \Omega})$
 is an $\Omega$-Lie-Yamaguti    algebra.
 We define a collection of linear maps as follows:
    \begin{small}
  \begin{align*}
  \left\{\begin{array}{cc}
 \varrho_{\alpha,\varsigma}:V\times L\rightarrow L,\rho_{{\alpha,\varsigma}}(u, x)=[T_\alpha u, x]+T_{\varsigma\alpha}(\rho(x)u+\Gamma_1(x,T_\alpha u) ), \\
\vartheta_{\alpha,\beta,\varsigma}: V \times V\times L\rightarrow L,\vartheta_{\alpha,\beta,\varsigma}(u,v,x)=\{x,T_\alpha u, T_\beta v\}- T_{\varsigma\alpha\beta}\big(D(x,T_\alpha u)v-\theta(x,T_\beta v)u+\Gamma_2(x,T_\alpha u, T_\beta v)\big)
 \end{array}\right\}_{\alpha,\beta,\varsigma\in\Omega}
 \end{align*}
  \end{small}

   \begin{theorem} \label{theorem:O-representation}
With the above notations,
 $(L,\{\varrho_{\alpha,\varsigma},\vartheta_{\alpha,\beta,\varsigma}\}_{\alpha,\beta,\varsigma\in \Omega})$ is a representation  of the $\Omega$-Lie-Yamaguti    algebra $(V,\{[\cdot, \cdot]_{\alpha,\beta}, \{\cdot, \cdot, \cdot\}_{\alpha,\beta,\gamma} \}_{\alpha,\beta,\gamma\in \Omega})$
 introduced in   Corollary  \ref{coro:O-LY}.
 \end{theorem}

    \begin{proof}
For any ~$x\in L, u,v,w\in V$ and $  \alpha,\beta,\gamma,\varsigma \in \Omega$,   we have
   \begin{align*}
&\mathcal{D}_{\alpha,\beta,\varsigma}(u,v,x)\\
=&\vartheta_{\beta,\alpha,\varsigma}(v,u,x)-\vartheta_{\alpha,\beta,\varsigma}(u,v,x)-\varrho_{\alpha\beta,\varsigma}([u,v]_{\alpha,\beta},x)+\varrho_{\alpha,\beta\varsigma}(u,\varrho_{\beta,\varsigma}(v,x))-
\varrho_{\beta,\alpha\varsigma}(v,\varrho_{\alpha,\varsigma}(u,x))\\
=&\{x,T_\beta v, T_\alpha u\}- T_{\varsigma\beta\alpha}\big(D(x,T_\beta v)u-\theta(x,T_\alpha u)v+\Gamma_2(x,T_\beta v, T_\alpha u)\big)-\\
  &\{x,T_\alpha u, T_\beta v\}+T_{\varsigma\alpha\beta}\big(D(x,T_\alpha u)v-\theta(x,T_\beta v)u+\Gamma_2(x,T_\alpha u, T_\beta v)\big)-\\
  &[ [T_{\alpha}u,T_{\beta}v], x]+T_{\varsigma\alpha\beta}\big(\rho(x)(\rho(T_\alpha u) v-\rho(T_\beta v) u+ \Gamma_1(T_\alpha u, T_\beta v))+\Gamma_1 (x,[T_{\alpha}u,T_{\beta}v])\big)+\\
  &[T_\alpha u, [T_\beta v, x]]+[T_\alpha u,T_{\varsigma\beta}(\rho(x)v+\Gamma_1\big(x,T_\beta v)\big)]+T_{\varsigma\beta\alpha}\Big(\rho\big([T_\beta v, x]+T_{\varsigma\beta}(\rho(x)v+\Gamma_1\big(x,T_\beta v))\big)u+\\
  &\Gamma_1\big([T_\beta v, x]+T_{\varsigma\beta}(\rho(x)v+\Gamma_1\big(x,T_\beta v)\big),T_\alpha u\big)\Big)-[T_\beta v, [T_\alpha u, x]]-[T_\beta v,T_{\varsigma\alpha}\big(\rho(x)u+\Gamma_1(x,T_\alpha u)\big)]-\\
  &T_{\varsigma\alpha\beta}\Big(\rho\big([T_\alpha u, x]+T_{\varsigma\alpha}(\rho(x)u+\Gamma_1\big(x,T_\alpha u))\big)v+\Gamma_1\big([T_\alpha u, x]+T_{\varsigma\alpha}(\rho(x)u+\Gamma_1\big(x,T_\alpha u)),T_\beta v\big)\Big)\\
  =&\{T_\alpha u, T_\beta v,x\}- T_{\alpha\beta\varsigma}\big(\theta(T_\beta v,x)u-\theta(T_\alpha u,x)v+\Gamma_2(T_\alpha u, T_\beta v,x)\big).
  \end{align*}
 It follows that
     \begin{align*}
&\vartheta_{\alpha\beta,\gamma,\varsigma}([u,v]_{\alpha,\beta},w,x)-\vartheta_{\alpha,\gamma,\beta\varsigma}(u,w,\varrho_{\beta,\varsigma}(v,x))+\vartheta_{\beta,\gamma,\alpha\varsigma}(v,w,\varrho_{\alpha,\varsigma}(u,x))\\
=&\{x, [T_\alpha u,T_\beta v], T_\gamma w\}- T_{\varsigma\alpha\beta\gamma}\Big(D(x,  [T_\alpha u,T_\beta v])w-\theta(x,T_\gamma w)\big(\rho(T_\alpha u) v-\rho(T_\beta v) u+ \Gamma_1(T_\alpha u, T_\beta v)\big)+\\
&\Gamma_2(x, [T_\alpha u,T_\beta v], T_\gamma w)\Big)-\{[T_\beta v, x]+T_{\varsigma\beta}(\rho(x)v+\Gamma_1\big(x,T_\beta v)\big),T_\alpha u, T_\gamma  w\}+ \\
&T_{\beta\varsigma\alpha\gamma}\Big(D\big([T_\beta v, x]+T_{\varsigma\beta}(\rho(x)v+\Gamma_1(x,T_\beta v)),T_\alpha u\big)w-\theta([T_\beta v, x]+T_{\varsigma\beta}(\rho(x)v+\Gamma_1\big(x,T_\beta v)\big),T_\gamma w)u+\\
&\Gamma_2([T_\beta v, x]+T_{\varsigma\beta}(\rho(x)v+\Gamma_1\big(x,T_\beta v)\big),T_\alpha u, T_\gamma w)\Big)+\{[T_\alpha u, x]+T_{\varsigma\alpha}(\rho(x)u+\Gamma_1(x,T_\alpha u)),T_\beta v, T_\gamma w\}- \\
&T_{\alpha\varsigma\beta\gamma}\Big(D\big([T_\alpha u, x]+T_{\varsigma\alpha}(\rho(x)u+\Gamma_1(x,T_\alpha u)),T_\beta v\big)w-\theta\big([T_\alpha u, x]+T_{\varsigma\alpha}(\rho(x)u+\Gamma_1(x,T_\alpha u)),T_\gamma w\big)v+\\
&\Gamma_2\big([T_\alpha u, x]+T_{\varsigma\alpha}(\rho(x)u+\Gamma_1(x,T_\alpha u)),T_\beta v, T_\gamma w\big)\Big)\\
=&0,\\
&\mathcal{D}_{\alpha,\beta,\gamma\varsigma}(u,v,\varrho_{\gamma,\varsigma}(w,x))-\varrho_{\gamma,\alpha\beta\varsigma}(w,\mathcal{D}_{\alpha,\beta,\varsigma}(u,v,x))-
\varrho_{\alpha\beta\gamma,\varsigma}(\{u,v,w\}_{\alpha,\beta,\gamma},x)\\
=&\big\{T_\alpha u, T_\beta v,[T_\gamma w, x]+T_{\varsigma\gamma}(\rho(x)w+\Gamma_1(x,T_\gamma w))\big\}- T_{\alpha\beta\gamma\varsigma}\Big(\theta\big(T_\beta v,[T_\gamma w, x]+T_{\varsigma\gamma}(\rho(x)w+\Gamma_1(x,T_\gamma w))\big)u-\\
&\theta\big(T_\alpha u,[T_\gamma w, x]+T_{\varsigma\gamma}(\rho(x)w+\Gamma_1(x,T_\gamma w))\big)v+\Gamma_2\big(T_\alpha u, T_\beta v,[T_\gamma w, x]+T_{\varsigma\gamma}(\rho(x)w+\Gamma_1(x,T_\gamma w))\big)\Big)-\\
&\big[T_\gamma w, \{T_\alpha u, T_\beta v,x\}- T_{\alpha\beta\varsigma}\big(\theta(T_\beta v,x)u-\theta(T_\alpha u,x)v+\Gamma_2(T_\alpha u, T_\beta v,x)\big)\big]-\\
&T_{\alpha\beta\varsigma\gamma}\Big(\rho\big(\{T_\alpha u, T_\beta v,x\}- T_{\alpha\beta\varsigma}\big(\theta(T_\beta v,x)u-\theta(T_\alpha u,x)v+\Gamma_2(T_\alpha u, T_\beta v,x)\big)\big)w+\Gamma_1\big(\{T_\alpha u, T_\beta v,x\}- \\
&T_{\alpha\beta\varsigma}(\theta(T_\beta v,x)u-\theta(T_\alpha u,x)v+\Gamma_2(T_\alpha u, T_\beta v,x)),T_\gamma w\big)\Big)-\big[\{T_\alpha u,T_\beta v,T_\gamma w\}, x\big]-\\
&T_{\varsigma\alpha\beta\gamma}\Big(\rho(x)\big(D(T_\alpha u, T_\beta v)w+\theta(T_\beta v, T_\gamma w)u-\theta(T_\alpha u, T_\gamma w)v+\Gamma_2(T_\alpha u, T_\beta v, T_\gamma w)\big)+\Gamma_1\big(x,\{T_\alpha u,T_\beta v,T_\gamma w\}) \Big)\\
=&0,\\
&\vartheta_{\alpha,\beta\gamma,\varsigma}(u,[v,w]_{\beta,\gamma},x)-\varrho_{\beta,\alpha\gamma\varsigma}(v,\vartheta_{\alpha,\gamma,\varsigma}(u,w,x))+\varrho_{\gamma,\alpha\beta\varsigma}(w,\vartheta_{\alpha,\beta,\varsigma}(u,v,x))\\
=&\{x,T_\alpha u, [T_\beta v,T_\gamma w]\}- T_{\varsigma\alpha\beta\gamma}\Big(D(x,T_\alpha u)\big(\rho(T_\beta v) w-\rho(T_\gamma w) v+ \Gamma_1(T_\beta v, T_\gamma w)\big)-\theta(x, [T_\beta v,T_\gamma w])u+\\
&\Gamma_2(x,T_\alpha u, [T_\beta v,T_\gamma w])\Big)-\big[T_\beta v, \{x,T_\alpha u, T_\gamma w\}- T_{\varsigma\alpha\gamma}\big(D(x,T_\alpha u)w-\theta(x,T_\gamma w)u+\Gamma_2(x,T_\alpha u, T_\gamma w)\big)\big]-\\
&T_{\varsigma\alpha\gamma\beta}\Big(\rho\big(\{x,T_\alpha u, T_\gamma w\}- T_{\varsigma\alpha\gamma}\big(D(x,T_\alpha u)w-\theta(x,T_\gamma w)u+\Gamma_2(x,T_\alpha u, T_\gamma w)\big)\big)v+\\
&\Gamma_1\big(\{x,T_\alpha u, T_\gamma w\}- T_{\varsigma\alpha\gamma}\big(D(x,T_\alpha u)w-\theta(x,T_\gamma w)u+\Gamma_2(x,T_\alpha u, T_\gamma w)\big),T_\beta v\big) \Big)+\\
&\big[T_\gamma w, \{x,T_\alpha u, T_\beta v\}- T_{\varsigma\alpha\beta}\big(D(x,T_\alpha u)v-\theta(x,T_\beta v)u+\Gamma_2(x,T_\alpha u, T_\beta v)\big)\big]+\\
&T_{\varsigma\alpha\beta\gamma}\Big(\rho\big(\{x,T_\alpha u, T_\beta v\}- T_{\varsigma\alpha\beta}\big(D(x,T_\alpha u)v-\theta(x,T_\beta v)u+\Gamma_2(x,T_\alpha u, T_\beta v)\big)\big)w+\\
&\Gamma_1(\{x,T_\alpha u, T_\beta v\}- T_{\varsigma\alpha\beta}\big(D(x,T_\alpha u)v-\theta(x,T_\beta v)u+\Gamma_2(x,T_\alpha u, T_\beta v)\big),T_\gamma w) \Big)\\
=&0,
  \end{align*}
Similarly,  for any ~$x\in L, u,v,w,o\in V$ and $  \alpha,\beta,\gamma,\varsigma,\tau \in \Omega$,  we have
     \begin{align*}
&\mathcal{D}_{\tau,\alpha,\beta\gamma\varsigma}(o,u,\vartheta_{\beta,\gamma,\varsigma}(v,w,x))-\vartheta_{\beta,\gamma,\tau\alpha\varsigma}(v,w,\mathcal{D}_{\tau,\alpha,\varsigma}(o,u,x))-\vartheta_{\tau\alpha\beta,\gamma,\varsigma}(\{o,u,v\}_{\varsigma,\alpha,\beta},w,x)-\\
&\vartheta_{\beta,\tau\alpha\gamma,\varsigma}(v,\{o,u,w\}_{\tau,\alpha,\gamma},x)=0,\\
&\vartheta_{\tau,\alpha\beta\gamma,\varsigma}(o,\{u,v,w\}_{\alpha,\beta,\gamma},x)-\vartheta_{\beta,\gamma,\tau\alpha\varsigma}(v,w,\vartheta_{\tau,\alpha,\varsigma}(o,u,x))+
\vartheta_{\alpha,\gamma,\tau\beta\varsigma}(u,w,\vartheta_{\tau,\beta,\varsigma}(o,v,x))-\\
&\mathcal{D}_{\alpha,\beta,\tau\gamma\varsigma}(u,v,\vartheta_{\tau,\gamma,\varsigma}(o,w,x))=0.
     \end{align*}
Therefore, we deduce that $(L,\{\varrho_{\alpha,\varsigma},\vartheta_{\alpha,\beta,\varsigma}\}_{\alpha,\beta,\varsigma\in \Omega})$ is a representation  of the $\Omega$-Lie-Yamaguti    algebra $(V,\{[\cdot, \cdot]_{\alpha,\beta}, \{\cdot, \cdot, \cdot\}_{\alpha,\beta,\gamma} \}_{\alpha,\beta,\gamma\in \Omega})$.
  \end{proof}

In the following of this section,  we assume that $\Omega$ is a commutative semigroup with unit $1\in \Omega$. The
unital condition of $\Omega$ is only useful in the coboundary operator of the cohomology at the
degree 0 level.

Let $(V; \rho, \theta)$ be a representation of a   Lie-Yamaguti   algebra $(L, [\cdot, \cdot], \{\cdot, \cdot, \cdot\})$,
 and $\{T_\alpha:V\rightarrow L\}_{\alpha\in\Omega}$ be an $\Gamma$-twisted    Rota-Baxter family.
  Then by Corollary \ref{coro:O-LY}, $\{T_\alpha:V\rightarrow L\}_{\alpha\in\Omega}$  induces an $\Omega$-Lie-Yamaguti    algebra $(V,\{[\cdot, \cdot]_{\alpha,\beta}, \{\cdot, \cdot, \cdot\}_{\alpha,\beta,\gamma} \}_{\alpha,\beta,\gamma\in \Omega})$.
Consequently, by Theorem \ref{theorem:O-representation},
we obtain a cochain complex of the $\Omega$-Lie-Yamaguti    algebra $(V,\{[\cdot, \cdot]_{\alpha,\beta}, \{\cdot, \cdot, \cdot\}_{\alpha,\beta,\gamma} \}_{\alpha,\beta,\gamma\in \Omega})$
 with the coefficients in the representation $(L,\{\varrho_{\alpha,\varsigma},\vartheta_{\alpha,\beta,\varsigma}\}_{\alpha,\beta,\varsigma\in \Omega})$:
$$\aligned
\xymatrix{
  C^{0}_{\Omega-\mathrm{LY}}(V, L)\ar[r]^-{\partial_{\Omega}} & C^{1}_{\Omega-\mathrm{LY}}(V, L)\ar[r]^-{\partial_{\Omega}} &  C^{(2,3)}_{\Omega-\mathrm{LY}}(V,L)\ar[r]^-{\partial_\Omega} \ar[d]^-{\partial^{*}_\Omega}& C^{(4,5)}_{\Omega-\mathrm{LY}}(V,L)\ar[r]^-{\partial_\Omega} & C^{(6,7)}_{\Omega-\mathrm{LY}}(V,L)\ar[r]^-{\partial_\Omega}&\cdots\\
 ~&~& C^{(3,4)}_{\Omega-\mathrm{LY}}(V,L) &  ~& ~&}
 \endaligned$$

%
Here, the  coboundary maps are defined as follows:

For $n\geq 1$,   $\partial_{\Omega}=(\partial_{\Omega,\mathrm{I}},\partial_{\Omega,\mathrm{II}}):C^{(2n,2n+1) }_{\mathrm{\Omega}-\mathrm{LY}}(V,L)\rightarrow C^{(2n+2,2n+3) }_{\mathrm{\Omega}-\mathrm{LY}}(V,L)$
is given by:
\begin{small}
\begin{align*}
&(\partial_{\Omega,\mathrm{I}}(f,g))_{\alpha_1,\ldots,\alpha_{2n+2}}(v_1,\ldots,v_{2n+2})\\
=&(-1)^n \Big([T_{\alpha_{2n+1}} v_{2n+1}, g_{\alpha_1,\ldots,\alpha_{2n},\alpha_{2n+2}}(v_1,\ldots,v_{2n},v_{2n+2})]+T_{\alpha_1\cdots\alpha_{2n}\alpha_{2n+2}\alpha_{2n+1}}\big(\rho\big(g_{\alpha_1,\ldots,\alpha_{2n},\alpha_{2n+2}}(v_1,\ldots,v_{2n},\\
&v_{2n+2})\big)v_{2n+1}+\Gamma_1\big(g_{\alpha_1,\ldots,\alpha_{2n},\alpha_{2n+2}}(v_1,\ldots,v_{2n},v_{2n+2}),T_{\alpha_{2n+1}} v_{2n+1}\big) \big)-\big[T_{\alpha_{2n+2}} v_{2n+2}, g_{\alpha_1,\ldots, \alpha_{2n+1}}(v_1,\ldots, v_{2n+1})\big]-\\
&T_{\alpha_1\cdots  \alpha_{2n+2} }\big(\rho\big(g_{\alpha_1,\ldots, \alpha_{2n+1}}(v_1,\ldots, v_{2n+1})\big)v_{2n+2}+\Gamma_1(g_{\alpha_1,\ldots, \alpha_{2n+1}}(v_1,\ldots, v_{2n+1}),T_{\alpha_{2n+2}} v_{2n+2}) \big)-\\
&g_{\alpha_1,\ldots, \alpha_{2n},\alpha_{2n+1}\alpha_{2n+2}}\big(v_1,\ldots, v_{2n}, \rho(T_{\alpha_{2n+1}} v_{2n+1})  v_{2n+2}-\rho(T_{\alpha_{2n+2}} v_{2n+2}) v_{2n+1}+ \Gamma_1(T_{\alpha_{2n+1}} v_{2n+1}, T_{\alpha_{2n+2}} v_{2n+2})
\big)\Big)+\\
&\sum_{k=1}^n(-1)^{k+1}\Big(\{T_{\alpha_{2k-1}} v_{2k-1}, T_{\alpha_{2k}} v_{2k},f_{\alpha_1,\ldots, \alpha_{2k-2},\alpha_{2k+1},\ldots, \alpha_{2n+2}}(v_1,\ldots,\widehat{v}_{2k-1},\widehat{v}_{2k}, \ldots, v_{2n+2})\}- \\
&T_{\alpha_{2k}\alpha_1\cdots\alpha_{2k-2}\alpha_{2k+1}\cdots\alpha_{2n+2}\alpha_{2k-1}}\big(\theta(T_{\alpha_{2k}} v_{2k},f_{\alpha_1,\ldots, \alpha_{2k-2},\alpha_{2k+1},\ldots, \alpha_{2n+2}}(v_1,\ldots,\widehat{v}_{2k-1},\widehat{v}_{2k}, \ldots, v_{2n+2}))v_{2k-1}-\\
&\theta(T_{\alpha_{2k-1}} v_{2k-1},f_{\alpha_1,\ldots, \alpha_{2k-2},\alpha_{2k+1},\ldots, \alpha_{2n+2}}(v_1,\ldots,\widehat{v}_{2k-1},\widehat{v}_{2k}, \ldots, v_{2n+2}))v_{2k}+\\
&\Gamma_2(T_{\alpha_{2k-1}} v_{2k-1}, T_{\alpha_{2k}} v_{2k},f_{\alpha_1,\ldots, \alpha_{2k-2},\alpha_{2k+1},\ldots, \alpha_{2n+2}}(v_1,\ldots,\widehat{v}_{2k-1},\widehat{v}_{2k}, \ldots, v_{2n+2}))\big)\Big)+\\
&\sum_{k=1}^n\sum_{j=2k+1}^{2n+2}(-1)^{k}f_{\alpha_1,\ldots, \alpha_{2k-2},\alpha_{2k+1},\ldots,\alpha_{2k-1}\alpha_{2k}\alpha_{j},\ldots, \alpha_{2n+2}}
\big(v_1,\ldots,\widehat{v}_{2k-1},\widehat{v}_{2k}, \ldots, D(T_{\alpha_{2k-1}} {v}_{2k-1}, T_{\alpha_{2k}} {v}_{2k})v_j+\\
&\theta(T_{\alpha_{2k}} {v}_{2k}, T_{\alpha_j} v_j){v}_{2k-1}-\theta(T_{\alpha_{2k-1}} {v}_{2k-1}, T_{\alpha_j} v_j){v}_{2k}+\Gamma_2(T_{\alpha_{2k-1}} {v}_{2k-1}, T_{\alpha_{2k}} {v}_{2k}, T_{\alpha_j} v_j), \ldots,v_{2n+2}\big);
\end{align*}
\begin{align*}
&(\partial_{\Omega,\mathrm{II}}(f,g))_{\alpha_1,\ldots,\alpha_{2n+3}}(v_1,\ldots,v_{2n+3})\\
=&(-1)^n \Big(\big\{g_{\alpha_1,\ldots, \alpha_{2n+1}}(v_1,\ldots, v_{2n+1}),T_{\alpha_{2n+2}} v_{2n+2}, T_{\alpha_{2n+3}} v_{2n+3}\big\}-\\
 &T_{\alpha_1\cdots\alpha_{2n+1}\alpha_{2n+2}\alpha_{2n+3}}\big(D(g_{\alpha_1,\ldots, \alpha_{2n+1}}(v_1,\ldots, v_{2n+1}),T_{\alpha_{2n+2}} v_{2n+2})v_{2n+3}-\\
 &\theta(g_{\alpha_1,\ldots, \alpha_{2n+1}}(v_1,\ldots, v_{2n+1}),T_{\alpha_{2n+3}} v_{2n+3})v_{2n+2}+\Gamma_2(g_{\alpha_1,\ldots, \alpha_{2n+1}}(v_1,\ldots, v_{2n+1}),T_{\alpha_{2n+2}} v_{2n+2}, T_{\alpha_{2n+3}} v_{2n+3})\big)-\\
 &\{g_{\alpha_1,\ldots, \alpha_{2n},\alpha_{2n+2}}(v_1,\ldots, v_{2n}, v_{2n+2})\big),T_{\alpha_{2n+1}} v_{2n+1}, T_{\alpha_{2n+3}} v_{2n+3}\}+ \\
 &T_{\alpha_1\cdots\alpha_{2n}\alpha_{2n+2}\alpha_{2n+1}\alpha_{2n+3}}\big(D(g_{\alpha_1,\ldots, \alpha_{2n},\alpha_{2n+2}}(v_1,\ldots, v_{2n}, v_{2n+2})\big),T_{\alpha_{2n+1}} v_{2n+1})v_{2n+3}-\\
 &\theta(g_{\alpha_1,\ldots, \alpha_{2n},\alpha_{2n+2}}(v_1,\ldots, v_{2n}, v_{2n+2})\big),T_{\alpha_{2n+3}} v_{2n+3})v_{2n+1}+\\
 &\Gamma_2(g_{\alpha_1,\ldots, \alpha_{2n},\alpha_{2n+2}}(v_1,\ldots, v_{2n}, v_{2n+2})\big),T_{\alpha_{2n+1}} v_{2n+1}, T_{\alpha_{2n+3}} v_{2n+3})\big)\Big)+\\
 &\sum_{k=1}^{n+1}(-1)^{k+1}\Big(\big\{T_{\alpha_{2k-1}} v_{2k-1}, T_{\alpha_{2k}} v_{2k},g_{\alpha_1,\ldots, \alpha_{2k-2},\alpha_{2k+1},\ldots, \alpha_{2n+3}}(v_1,\ldots,\widehat{v}_{2k-1},\widehat{v}_{2k}, \ldots, v_{2n+3})\big\}- \\
&T_{\alpha_{2k}\alpha_1\cdots\alpha_{2k-2}\alpha_{2k+1}\cdots\alpha_{2n+3}\alpha_{2k-1}}\big(\theta(T_{\alpha_{2k}} v_{2k},g_{\alpha_1,\ldots, \alpha_{2k-2},\alpha_{2k+1},\ldots, \alpha_{2n+3}}(v_1,\ldots,\widehat{v}_{2k-1},\widehat{v}_{2k}, \ldots, v_{2n+3}))v_{2k-1}-\\
&\theta(T_{\alpha_{2k-1}} v_{2k-1},g_{\alpha_1,\ldots, \alpha_{2k-2},\alpha_{2k+1},\ldots, \alpha_{2n+3}}(v_1,\ldots,\widehat{v}_{2k-1},\widehat{v}_{2k}, \ldots, v_{2n+3}))v_{2k}+\\
&\Gamma_2(T_{\alpha_{2k-1}} v_{2k-1}, T_{\alpha_{2k}} v_{2k},g_{\alpha_1,\ldots, \alpha_{2k-2},\alpha_{2k+1},\ldots, \alpha_{2n+3}}(v_1,\ldots,\widehat{v}_{2k-1},\widehat{v}_{2k}, \ldots, v_{2n+3})\big)\Big)+\\
&\sum_{k=1}^{n+1}\sum_{j=2k+1}^{2n+3}(-1)^{k}g_{\alpha_1,\ldots, \alpha_{2k-2},\alpha_{2k+1},\ldots,\alpha_{2k-1}\alpha_{2k}\alpha_{j},\ldots, \alpha_{2n+3}}
(v_1,\ldots,\widehat{v}_{2k-1},\widehat{v}_{2k}, \ldots, D(T_{\alpha_{2k-1}} {v}_{2k-1}, T_{\alpha_{2k}} {v}_{2k})v_j+\\
&\theta(T_{\alpha_{2k}} {v}_{2k}, T_{\alpha_j} v_j){v}_{2k-1}-\theta(T_{\alpha_{2k-1}} {v}_{2k-1}, T_{\alpha_j} v_j){v}_{2k}+\Gamma_2(T_{\alpha_{2k-1}} {v}_{2k-1}, T_{\alpha_{2k}} {v}_{2k}, T_{\alpha_j} v_j), \ldots,v_{2n+3})\big).
\end{align*}
\end{small}
The map $\partial^*_{\Omega}=(\partial^*_{\Omega,\mathrm{I}},\partial^*_{\Omega,\mathrm{II}}):C^{(2,3) }_{\mathrm{\Omega}-\mathrm{LY}}(V,L)\rightarrow C^{(3,4) }_{\mathrm{\Omega}-\mathrm{LY}}(V,L)$ is given by:
\begin{align*}
&(\partial^*_{\Omega,\mathrm{I}}(f,g))_{\alpha_1,\alpha_2,\alpha_3}(v_1,v_2,v_3)\\
=&-[T_{\alpha_1} v_1, f_{\alpha_2,\alpha_3}(v_2,v_3)]-T_{\alpha_2\alpha_3\alpha_1}\big(\rho(f_{\alpha_2,\alpha_3}(v_2,v_3))v_1+\Gamma_1(f_{\alpha_2,\alpha_3}(v_2,v_3),T_{\alpha_1} v_1) \big)-\\
&[T_{\alpha_2} v_2, f_{\alpha_3,\alpha_1}(v_3,v_1)]-T_{\alpha_3\alpha_1\alpha_2}\big(\rho(f_{\alpha_3,\alpha_1}(v_3,v_1))  v_2+\Gamma_1(f_{\alpha_3,\alpha_1}(v_3,v_1),T_{\alpha_2} v_2) \big)-\\
&[T_{\alpha_2} v_3, f_{\alpha_1,\alpha_2}(v_1,v_2)]-T_{\alpha_1\alpha_2\alpha_3}\big(\rho(f_{\alpha_1,\alpha_2}(v_1,v_2))  v_3+\Gamma_1(f_{\alpha_1,\alpha_2}(v_1,v_2),T_{\alpha_3} v_3) \big)+\\
%
&f_{\alpha_1\alpha_2,\alpha_3}\big(\rho(T_{\alpha_1} v_1) v_2-\rho(T_{\alpha_2} v_2) v_1+ \Gamma_1(T_{\alpha_1} v_1, T_{\alpha_2} v_2),v_3\big)+\\
&f_{\alpha_2\alpha_3,\alpha_1}(\rho(T_{\alpha_2} v_2) v_3-\rho(T_{\alpha_3} v_3) v_2+ \Gamma_1(T_{\alpha_2} v_2, T_{\alpha_3} v_3),v_1)+\\
&f_{\alpha_3\alpha_1,\alpha_2}(\rho(T_{\alpha_3} v_3) v_1-\rho(T_{\alpha_1} v_1) v_3+ \Gamma_1(T_{\alpha_3} v_3, T_{\alpha_1} v_1),v_2)+\\
&g_{\alpha_1,\alpha_2,\alpha_3}(v_1,v_2,v_3)+g_{\alpha_2,\alpha_3,\alpha_1}(v_2, v_3,v_1)+g_{\alpha_3,\alpha_1,\alpha_2}(v_3,v_1,v_2),\\
&(\partial^*_{\Omega,\mathrm{II}}(f,g))_{\alpha_1,\alpha_2,\alpha_3,\alpha_4}(v_1,v_2,v_3,v_4)\\
=&\{f_{\alpha_2,\alpha_3}(v_2,v_3),T_{\alpha_1}v_1, T_{\alpha_4} v_4\}- T_{\alpha_1\alpha_4,\alpha_2\alpha_3}\Big(D(f_{\alpha_2,\alpha_3}(v_2,v_3),T_{\alpha_1} v_1)v_4-\theta(f_{\alpha_2,\alpha_3}(v_2,v_3),T_{\alpha_4} v_4)v_1+\\
&\Gamma_2(f_{\alpha_2,\alpha_3}(v_2,v_3),T_{\alpha_1} v_1, T_{\alpha_4} v_4)\Big)+\{f_{\alpha_3,\alpha_1}(v_3,v_1),T_{\alpha_2}v_2, T_{\alpha_4} v_4\}- T_{\alpha_2\alpha_4,\alpha_3\alpha_1}\Big(D(f_{\alpha_3,\alpha_1}(v_3,v_1),T_{\alpha_2} v_2)v_4-\\
&\theta(f_{\alpha_3,\alpha_1}(v_3,v_1),T_{\alpha_4} v_4)v_2+\Gamma_2(f_{\alpha_3,\alpha_1}(v_3,v_1),T_{\alpha_2} v_2, T_{\alpha_4} v_4)\Big)+\{f_{\alpha_1,\alpha_2}(v_1,v_2),T_{\alpha_3}v_3, T_{\alpha_4} v_4\}- \\ &T_{\alpha_3\alpha_4,\alpha_1\alpha_2}\Big(D(f_{\alpha_1,\alpha_2}(v_1,v_2),T_{\alpha_3} v_3)v_4-\theta(f_{\alpha_1,\alpha_2}(v_1,v_2),T_{\alpha_4} v_4)v_3+\Gamma_2(f_{\alpha_1,\alpha_2}(v_1,v_2),T_{\alpha_3} v_3, T_{\alpha_4} v_4)\Big)+\\
&g_{\alpha_1\alpha_2,\alpha_3,\alpha_4}(\rho(T_{\alpha_1} v_1) v_2-\rho(T_{\alpha_2} v_2) v_1+ \Gamma_1(T_{\alpha_1} v_1, T_{\alpha_2} v_2),v_3,v_4)+\\
&g_{\alpha_2\alpha_3,\alpha_1,\alpha_4}(\rho(T_{\alpha_2} v_2) v_3-\rho(T_{\alpha_3} v_3) v_2+ \Gamma_1(T_{\alpha_2} v_2, T_{\alpha_3} v_3),v_1,v_4)+\\
&g_{\alpha_3\alpha_1,\alpha_2,\alpha_4}(\rho(T_{\alpha_3} v_3) v_1-\rho(T_{\alpha_1} v_1) v_3+ \Gamma_1(T_{\alpha_3} v_3, T_{\alpha_1} v_1),v_2,v_4).
\end{align*}

The map $\partial_{\Omega}=(\partial_{\Omega,\mathrm{I}},\partial_{\Omega,\mathrm{II}}):C^{1 }_{\mathrm{\Omega}-\mathrm{LY}}(V,L)\rightarrow C^{(2,3) }_{\mathrm{\Omega}-\mathrm{LY}}(V,L)$ is given by:
\begin{align*}
&(\partial_{\Omega,\mathrm{I}}f)_{\alpha_1,\alpha_2}(v_1,v_2)\\
=&[T_{\alpha_1} v_1, f_{\alpha_2}(v_2)]+T_{\alpha_1\alpha_2}\big(\rho(f_{\alpha_2}(v_2))v_1+\Gamma_1(f_{\alpha_2}(v_2),T_{\alpha_1} v_1) \big)-[T_{\alpha_2} v_2, f_{\alpha_1}(v_1)]-\\
&T_{\alpha_2\alpha_1}\big(\rho(f_{\alpha_1}(v_1))v_2+\Gamma_1(f_{\alpha_1}(v_1),T_{\alpha_2} v_2) \big)-f_{\alpha_1\alpha_2}\big(\rho(T_{\alpha_1} v_1) v_2-\rho(T_{\alpha_2} v_2) v_1+ \Gamma_1(T_{\alpha_1} v_1, T_{\alpha_2} v_2)\big),\\
&(\partial_{\Omega,\mathrm{II}}f)_{\alpha_1,\alpha_2,\alpha_3}(v_1,v_2,v_3)\\
=&\{T_{\alpha_1} v_1, T_{\alpha_2} v_2,f_{\alpha_3}(v_3)\}- T_{\alpha_1\alpha_2\alpha_3}\big(\theta(T_{\alpha_2} v_2,f_{\alpha_3}(v_3))v_1-\theta(T_{\alpha_1} v_1,f_{\alpha_3}(v_3))v_2+\Gamma_2(T_{\alpha_1} v_1, T_{\alpha_2} v_2,f_{\alpha_3}(v_3))\big)+\\
&\{f_{\alpha_1}(v_1),T_{\alpha_2} v_2, T_{\alpha_3} v_3\}- T_{\alpha_1\alpha_2\alpha_3}\big(D(f_{\alpha_1}(v_1),T_{\alpha_2}v_2)v_3-\theta(f_{\alpha_1}(v_1),T_{\alpha_3} v_3)v_2+\Gamma_2(f_{\alpha_1}(v_1),T_{\alpha_2} v_2, T_{\alpha_3} v_3)\big)-\\
&\{f_{\alpha_2}(v_2),T_{\alpha_1} v_1, T_{\alpha_3} v_3\}+ T_{\alpha_2\alpha_1\alpha_3}\big(D(f_{\alpha_2}(v_2),T_{\alpha_1}v_1)v_3-\theta(f_{\alpha_2}(v_2),T_{\alpha_3} v_3)v_1+\Gamma_2(f_{\alpha_2}(v_2),T_{\alpha_1} v_1, T_{\alpha_3} v_3)\big)-\\
&f_{\alpha_1\alpha_2\alpha_3}\big(D(T_{\alpha_1} v_1, T_{\alpha_2}  v_2)v_3+\theta(T_{\alpha_2}  v_2, T_{\alpha_3}  v_3)v_1-\theta(T_{\alpha_1}  v_1, T_{\alpha_3}  v_3)v_2+\Gamma_2(T_{\alpha_1}  v_1, T_{\alpha_2}  v_2, T_{\alpha_3}  v_3)\big).
%
%
\end{align*}
The map $\partial_{\Omega}: C^{0}_{\mathrm{\Omega}-\mathrm{LY}}(V,L)\rightarrow C^{1 }_{\mathrm{\Omega}-\mathrm{LY}}(V,L)$ is given by:
\begin{align}
(\partial_{\Omega}(a,b))_{\alpha}u:=T_{\alpha}(D(a,b)u+\Gamma_2(a,b,T_{\alpha} u))-\{a,b,T_{\alpha} u\},  \label{6.1}
\end{align}
 for all $(a,b)\in C^{0}_{\mathrm{\Omega}-\mathrm{LY}}(V,L):=\wedge ^2 L ~\text{and}~ u\in V.$\\

The above  cochain complex  at the degree 0 level  is justified by the following lemma.

 \begin{lemma} \label{prop:1-cocycle}
 With the notations above,
$\partial_{\Omega}(\partial_{\Omega}(a,b))=0$,  that is the composition  $
C^{0}_{\Omega-\mathrm{LY}}(V, L)\stackrel{\partial_{\Omega}}{\longrightarrow} C^{1}_{\Omega-\mathrm{LY}}(V, L)\stackrel{\partial_{\Omega}}{\longrightarrow} C^{(2,3)}_{\Omega-\mathrm{LY}}(V, L)$ is the zero map.
 \end{lemma}
  \begin{proof}
  For any $v_1,v_2\in V$ and $\alpha_1,\alpha_2\in \Omega$, by Eqs. \eqref{2.3}, \eqref{2.6}, \eqref{2.16} and \eqref{3.1}, we have
   \begin{small}
  \begin{align*}
&(\partial_{\Omega,\mathrm{I}}(\partial_{\Omega}(a,b)))_{\alpha_1,\alpha_2}(v_1,v_2)\\
=&[T_{\alpha_1} v_1, (\partial_{\Omega}(a,b))_{\alpha_2}(v_2)]+T_{\alpha_1\alpha_2}\big(\rho((\partial_{\Omega}(a,b))_{\alpha_2}(v_2))v_1+\Gamma_1((\partial_{\Omega}(a,b))_{\alpha_2}(v_2),T_{\alpha_1} v_1) \big)-\\
&\big[T_{\alpha_2} v_2, (\partial_{\Omega}(a,b))_{\alpha_1}(v_1)\big]-T_{\alpha_2\alpha_1}\big(\rho((\partial_{\Omega}(a,b))_{\alpha_1}(v_1))v_2+\Gamma_1((\partial_{\Omega}(a,b))_{\alpha_1}(v_1),T_{\alpha_2} v_2) \big)-\\
&(\partial_{\Omega}(a,b))_{\alpha_1\alpha_2}\big(\rho(T_{\alpha_1} v_1) v_2-\rho(T_{\alpha_2} v_2) v_1+ \Gamma_1(T_{\alpha_1} v_1, T_{\alpha_2} v_2)\big)\\
=&\big[T_{\alpha_1} v_1, T_{\alpha_2} D(a,b)v_2+T_{\alpha_2} \Gamma_2(a,b,T_{\alpha_2} v_2) -\{a,b,T_{\alpha_2} v_2\}\big]+T_{\alpha_1\alpha_2}\Big(\rho\big(T_{\alpha_2} D(a,b)v_2+T_{\alpha_2} \Gamma_2(a,b,T_{\alpha_2} v_2) -\\
&\{a,b,T_{\alpha_2} v_2\}\big)v_1+\Gamma_1(T_{\alpha_2} D(a,b)v_2+T_{\alpha_2} \Gamma_2(a,b,T_{\alpha_2} v_2) -\{a,b,T_{\alpha_2} v_2\},T_{\alpha_1} v_1)\Big)-\\
&[T_{\alpha_2} v_2, T_{\alpha_1} D(a,b)v_1+T_{\alpha_1} \Gamma_2(a,b,T_{\alpha_1} v_1) -\{a,b,T_{\alpha_1} v_1\}]-T_{\alpha_1\alpha_2}\Big(\rho\big(T_{\alpha_1} D(a,b)v_1+T_{\alpha_1} \Gamma_2(a,b,T_{\alpha_1} v_1) -\\
&\{a,b,T_{\alpha_1} v_1\}\big)v_2+\Gamma_1(T_{\alpha_1} D(a,b)v_1+T_{\alpha_1} \Gamma_2(a,b,T_{\alpha_1} v_1) -\{a,b,T_{\alpha_1} v_1\},T_{\alpha_2} v_2) \Big)-T_{\alpha_1\alpha_2}\Big(D(a,b)\big(\rho(T_{\alpha_1} v_1) v_2-\\
&\rho(T_{\alpha_2} v_2) v_1+ \Gamma_1(T_{\alpha_1} v_1, T_{\alpha_2} v_2)\big)+\Gamma_2(a,b, [T_{\alpha_1}v_1,T_{\alpha_2}v_2])\Big)+\{a,b, [T_{\alpha_1}v_1,T_{\alpha_2}v_2]\}\\
=& [T_{\alpha_1} v_1, T_{\alpha_2} D(a,b)v_2]+ [T_{\alpha_1} v_1,T_{\alpha_2} \Gamma_2(a,b,T_{\alpha_2} v_2)  ]-[T_{\alpha_2} v_2, T_{\alpha_1} D(a,b)v_1]-[T_{\alpha_2} v_2,T_{\alpha_1} \Gamma_2(a,b,T_{\alpha_1} v_1)]+\\
&T_{\alpha_1\alpha_2}\Big(\rho\big(T_{\alpha_2} D(a,b)v_2\big)v_1+\rho\big(T_{\alpha_2} \Gamma_2(a,b,T_{\alpha_2} v_2)\big)v_1 -\rho\big(\{a,b,T_{\alpha_2} v_2\}\big)v_1+\Gamma_1(T_{\alpha_2} D(a,b)v_2,T_{\alpha_1} v_1)+\\
&\Gamma_1(T_{\alpha_2} \Gamma_2(a,b,T_{\alpha_2} v_2),T_{\alpha_1} v_1) -\Gamma_1(\{a,b,T_{\alpha_2} v_2\},T_{\alpha_1} v_1) - \rho\big(T_{\alpha_1} D(a,b)v_1\big)v_2-\rho\big(T_{\alpha_1} \Gamma_2(a,b,T_{\alpha_1} v_1)\big)v_2 +\\
&\rho\big(\{a,b,T_{\alpha_1} v_1\}\big)v_2-\Gamma_1\big(T_{\alpha_1} D(a,b)v_1,T_{\alpha_2} v_2\big)-\Gamma_1\big(T_{\alpha_1} \Gamma_2(a,b,T_{\alpha_1} v_1),T_{\alpha_2} v_2\big) +\Gamma_1\big(\{a,b,T_{\alpha_1} v_1\},T_{\alpha_2} v_2\big)-\\
&D(a,b)\big(\rho(T_{\alpha_1} v_1) v_2\big)+D(a,b)\big(\rho(T_{\alpha_2} v_2) v_1\big)-D(a,b)\big( \Gamma_1(T_{\alpha_1} v_1, T_{\alpha_2} v_2)\big)-\Gamma_2(a,b, [T_{\alpha_1}v_1,T_{\alpha_2}v_2])\Big)\\
=& [T_{\alpha_1} v_1, T_{\alpha_2} D(a,b)v_2]+ [T_{\alpha_1} v_1,T_{\alpha_2} \Gamma_2(a,b,T_{\alpha_2} v_2)  ]-[T_{\alpha_2} v_2, T_{\alpha_1} D(a,b)v_1]-[T_{\alpha_2} v_2,T_{\alpha_1} \Gamma_2(a,b,T_{\alpha_1} v_1)]+\\
&T_{\alpha_1\alpha_2}\Big(\rho\big(T_{\alpha_2} D(a,b)v_2\big)v_1+\rho\big(T_{\alpha_2} \Gamma_2(a,b,T_{\alpha_2} v_2)\big)v_1  +\Gamma_1(T_{\alpha_2} D(a,b)v_2,T_{\alpha_1} v_1)+\Gamma_1(T_{\alpha_2} \Gamma_2(a,b,T_{\alpha_2} v_2),T_{\alpha_1} v_1)  -\\
& \rho\big(T_{\alpha_1} D(a,b)v_1\big)v_2-\rho\big(T_{\alpha_1} \Gamma_2(a,b,T_{\alpha_1} v_1)\big)v_2 -\Gamma_1\big(T_{\alpha_1} D(a,b)v_1,T_{\alpha_2} v_2\big)-\Gamma_1\big(T_{\alpha_1} \Gamma_2(a,b,T_{\alpha_1} v_1),T_{\alpha_2} v_2\big) -\\
&\rho(T_{\alpha_1} v_1)\big(D(a,b) v_2\big)+\rho(T_{\alpha_2} v_2)\big(D(a,b)v_1\big)-\rho(T_{\alpha_1} v_1)  \Gamma_2(a,b, T_{\alpha_2} v_2) +\rho(T_{\alpha_2} v_2)\Gamma_2(a,b,  T_{\alpha_1}v_1)\Big)\\
=&0.
\end{align*}
\end{small}
Similarly, we can show that $(\partial_{\Omega,\mathrm{II}}(\partial_{\Omega}(a,b)))_{\alpha_1,\alpha_2,\alpha_3}(v_1,v_2,v_3)=0$ for all $v_1,v_2,v_3\in V$ and $\alpha_1,\alpha_2,\alpha_3\in \Omega$.
 \end{proof}

\begin{defn}
Let
\begin{align*}
Z^{(2n,2n+1) }_{\mathrm{\Omega}-\mathrm{LY}}(V,L)=&\text{span}\{(f,g)\in C^{(2n,2n+1) }_{\mathrm{\Omega}-\mathrm{LY}}(V,L)~|~\partial_{\Omega}(f,g)=0\}, n\geq 2;\\
Z^{(2,3) }_{\mathrm{\Omega}-\mathrm{LY}}(V,L)=&\text{span}\{(f,g)\in C^{(2,3) }_{\mathrm{\Omega}-\mathrm{LY}}(V,L)~|~\partial_{\Omega}(f,g)=0=\partial^*_{\Omega}(f,g)\} \ ~\ ~  \text{and}\\
Z^{1 }_{\mathrm{\Omega}-\mathrm{LY}}(V,L)=&\text{span}\{f\in C^{1 }_{\mathrm{\Omega}-\mathrm{LY}}(V,L)~|~\partial_{\Omega}(f)=0\}.
 \end{align*}
 An element of $Z^{(2n,2n+1) }_{\mathrm{\Omega}-\mathrm{LY}}(V,L)$ is called a $(2n,2n+1)$-cocycle, and an element of $Z^{1}_{\mathrm{\Omega}-\mathrm{LY}}(V,L)$  is referred to a $1$-cocycle.
Let
\begin{align*}
&B^{(2n,2n+1) }_{\mathrm{\Omega}-\mathrm{LY}}(V,L)=\partial_{\Omega}C^{(2n-2,2n-1) }_{\mathrm{\Omega}-\mathrm{LY}}(V,L),\ ~~n\geq2; \\
&B^{(2,3) }_{\mathrm{\Omega}-\mathrm{LY}}(V,L)=\partial_{\Omega}C^{1 }_{\mathrm{\Omega}-\mathrm{LY}}(V,L) ~\text{and}~ B^{1 }_{\mathrm{\Omega}-\mathrm{LY}}(V,L)=\partial_{\Omega}C^{0 }_{\mathrm{\Omega}-\mathrm{LY}}(V,L).
 \end{align*}
An element of $B^{(2n,2n+1) }_{\mathrm{\Omega}-\mathrm{LY}}(V,L)$ is called a $(2n,2n+1)$-coboundary, while an element of $B^{1 }_{\mathrm{\Omega}-\mathrm{LY}}(V,L)$ is called a $1$-coboundary.
The corresponding cohomology groups,
$$H^{1 }_{\mathrm{\Omega}-\mathrm{LY}}(V,L)=\frac{Z^{1 }_{\mathrm{\Omega}-\mathrm{LY}}(V,L)}{B^{1}_{\mathrm{\Omega}-\mathrm{LY}}(V,L)} ~\text{and} ~H^{(2n,2n+1) }_{\mathrm{\Omega}-\mathrm{LY}}(V,L)=\frac{Z^{(2n,2n+1) }_{\mathrm{\Omega}-\mathrm{LY}}(V,L)}{B^{(2n,2n+1) }_{\mathrm{\Omega}-\mathrm{LY}}(V,L)}, \text{where}~ n\geq 1, $$
are referred to as the cohomology of the $\Gamma$-twisted    Rota-Baxter family  $\{T_\alpha\}_{\alpha\in\Omega}$.
\end{defn}

   \begin{remark}
 (i)  Let $\{T_\alpha:V\rightarrow L\}_{\alpha\in\Omega}$  be a relative  Rota-Baxter family  (see Example \ref{exam:2.5}). Then $(V,\big\{[\cdot, \cdot]_{\alpha,\beta}, \{\cdot, \cdot, \cdot\}_{\alpha,\beta,\gamma} \big\}_{\alpha,\beta,\gamma\in \Omega})$ is
  an $\Omega$-Lie-Yamaguti    algebra, where
 \begin{align*}
[u, v]_{\alpha,\beta}=& \rho(T_\alpha u) v-\rho(T_\beta v) u,\\
\{u,v,w\}_{\alpha,\beta,\gamma}=& D(T_\alpha u, T_\beta v)w+\theta(T_\beta v, T_\gamma w)u-\theta(T_\alpha u, T_\gamma w)v,
  \end{align*}
  for all $u,v,w\in V$ and $\alpha,\beta,\gamma\in \Omega.$

Moreover,  $(L,\big\{\rho_{\alpha,\varsigma},\theta_{\alpha,\beta,\varsigma}\big\}_{\alpha,\beta,\varsigma\in \Omega})$ is a representation  of the $\Omega$-Lie-Yamaguti    algebra $(V,\big\{[\cdot, \cdot]_{\alpha,\beta},$ $ \{\cdot, \cdot, \cdot\}_{\alpha,\beta,\gamma} \big\}_{\alpha,\beta,\gamma\in \Omega})$, where,
\begin{small}
  \begin{align*}
  \left\{\begin{array}{cc}
 \rho_{\alpha,\varsigma}:V\times L\rightarrow L,\rho_{{\alpha,\varsigma}}(u, x)=[T_\alpha u, x]+T_{\varsigma\alpha}(\rho(x)u), \\
\theta_{\alpha,\beta,\varsigma}: V \times V\times L\rightarrow L,\vartheta_{\alpha,\beta,\varsigma}(u,v,x)=\{x,T_\alpha u, T_\beta v\}- T_{\varsigma\alpha\beta}\big(D(x,T_\alpha u)v-\theta(x,T_\beta v)u\big)
 \end{array}\right\}_{\alpha,\beta,\varsigma\in\Omega}.
 \end{align*}
  \end{small}
 The corresponding cohomology is called the cohomology of the relative  Rota-Baxter family  $\{T_\alpha:V\rightarrow L\}_{\alpha\in\Omega}$.

 (ii) The cohomology of the $\Gamma$-twisted  Rota-Baxter family  defined as above is the
cohomology of the generalized Reynolds operator (see Example \ref{exam:generalized Reynolds operator})  when the semigroup $\Omega$ is trivial.
 \end{remark}

At the end of this section, we study the    deformations of $\Gamma$-twisted    Rota-Baxter families.
We show that the cohomology of  $\Gamma$-twisted    Rota-Baxter families introduced above govern such    deformations.

 Let $(V; \rho, \theta)$ be a representation of a   Lie-Yamaguti   algebra $(L, [\cdot, \cdot], \{\cdot, \cdot, \cdot\})$,
 and $T=\{T_\alpha:V\rightarrow L\}_{\alpha\in\Omega}$ be an $\Gamma$-twisted    Rota-Baxter family.

 A linear     deformation of $T$ consists of a sum $T^t=T+tT^1$, where $T^1=\{T^1_\alpha:V\rightarrow L\}_{\alpha\in\Omega}$, such that $T^t$ is an $\Gamma$-twisted Rota-Baxter family, for all $t$.
  That is
\begin{small}
\begin{align*}
&[T^t_\alpha u, T^t_\beta v]=T^t_{\alpha\beta}\big(\rho(T^t_\alpha u) v-\rho(T^t_\beta v)u+\Gamma_1(T^t_\alpha u, T^t_\beta v)\big), \\
&\{T^t_\alpha u, T^t_\beta v, T^t_\gamma w\}=T^t_{\alpha\beta\gamma}\big(D(T^t_\alpha u, T^t_\beta v)w-\theta(T^t_\alpha u, T^t_\gamma w)v+\theta(T^t_\beta v, T^t_\gamma w)u+\Gamma_2(T^t_\alpha u, T^t_\beta v,T^t_\gamma w)\big),
\end{align*}
\end{small}
for all $u,v,w\in V$ and $\alpha,\beta,\gamma\in \Omega.$
By expanding the above identity and comparing coefficients
of various powers of $t^1$, we have
\begin{small}
\begin{align}
&[T^1_\alpha u, T_\beta v]+[T_\alpha u, T^1_\beta v]=T^1_{\alpha\beta}\big(\rho(T_\alpha u ) v-\rho(T_\beta v)u+\Gamma_1(T_\alpha u, T_\beta v)\big)+ \label{6.2} \\
& T_{\alpha\beta}\big(\rho(T^1_\alpha u ) v-\rho(T^1_\beta v)u+\Gamma_1(T^1_\alpha u, T_\beta v)+\Gamma_1(T_\alpha u, T^1_\beta v)\big), \nonumber\\
&\{T^1_\alpha u, T_\beta v, T_\gamma w\}+\{T_\alpha u, T^1_\beta v, T_\gamma w\}+\{T_\alpha u, T_\beta v, T^1_\gamma w\}\label{6.3} \\
&=T^1_{\alpha\beta\gamma}\big(D(T_\alpha u, T_\beta v)w-\theta(T_\alpha u, T_\gamma w)v+\theta(T_\beta v, T_\gamma w)u+\Gamma_2(T_\alpha u, T_\beta v,T_\gamma w)\big)+\nonumber\\
&  T_{\alpha\beta\gamma}\Big(D(T^1_\alpha u, T_\beta v)w+D(T_\alpha u, T^1_\beta v)w-\theta(T^1_\alpha u, T_\gamma w)v-\theta(T_\alpha u, T^1_\gamma w)v+\theta(T^1_\beta v, T_\gamma w)u+\nonumber\\
& \theta(T_\beta v, T^1_\gamma w)u+  \Gamma_2(T^1_\alpha u, T_\beta v,T_\gamma w)+\Gamma_2(T_\alpha u, T^1_\beta v,T_\gamma w)+\Gamma_2(T_\alpha u, T_\beta v,T^1_\gamma w)\Big).\nonumber
\end{align}
\end{small}
It follows from Eqs. \eqref{6.2} and  \eqref{6.3}     that   $(\partial_{\Omega,\mathrm{I}}T^1)_{\alpha,\beta}(u,v)=0$ and  $(\partial_{\Omega,\mathrm{II}}T^1)_{\alpha,\beta,\gamma}(u,v,w)=0$, for all  $u,v,w\in V$ and $\alpha,\beta,\gamma\in \Omega$.
 That is $T^1\in C^{1 }_{\mathrm{\Omega}-\mathrm{LY}}(V,L)$ is a 1-cocycle in the cochain complex of the $\Gamma$-twisted   Rota-Baxter family  $\{T_\alpha:V\rightarrow L\}_{\alpha\in\Omega}$.
The 1-cocycle $T^1$ is called the infinitesimal of the deformation $T^t$.

\begin{defn}
 Two linear   deformations $T^t=T+tT^1$ and $ \overline{T}^t=T+t \overline{T}^1$
 of $T$ are equivalent if there exists a pair $(a,b)\in  \wedge ^2 L$   such that the maps
 \begin{align*}
\big\{\zeta^t_\alpha= \mathrm{id}_V+t\big(D(a,b)(\cdot)+\Gamma_2(a,b,T_\alpha(\cdot))\big): V[t]\rightarrow V[t] \big\}_{\alpha\in \Omega},~ \ ~ \ ~ \ ~ \ ~ \eta^t=\mathrm{id}_L+t\{a,b,\cdot\}:L[t]\rightarrow L[t]
\end{align*}
  satisfy  $\eta^t\circ T_\alpha^t= \overline{T}_\alpha^t\circ \zeta^t_\alpha$ for all $\alpha\in \Omega$. Here the map $\eta^t$ doesn't depend on $\alpha$.
\end{defn}

By equating coefficients of $t$ in both sides of $\eta^t( T_\alpha^tu)= \overline{T}_\alpha^t\zeta^t_\alpha(u)$ for all $u\in V$, we have
\begin{align*}
&T_\alpha^1 u+\{a,b,T_\alpha u\}=\overline{T}_\alpha^1 u+T_\alpha\big(D(a,b)u+\Gamma_2(a,b,T_\alpha u)\big),
\end{align*}
that is
\begin{align*}
&T_\alpha^1 u-\overline{T}_\alpha^1 u=T_\alpha\big(D(a,b)u+\Gamma_2(a,b,T_\alpha u)\big)-\{a,b,T_\alpha u\}=(\partial_{\Omega}(a,b))_{\alpha}u.
\end{align*}
As a summary, we get the following theorem.
\begin{theorem}
  Let $T^t=T+tT^1$ and $ \overline{T}^t=T+t \overline{T}^1$ be two equivalent linear  deformations of the  $\Gamma$-twisted    Rota-Baxter family $T=\{T_\alpha\}_{\alpha\in\Omega}$,
    then their infinitesimals $T^1$
and $\overline{T}^1$ are in the same cohomology class.
\end{theorem}

Regarding the formal deformation of $\Gamma$-twisted   Rota-Baxter families, analogous results are obtained.
 \begin{remark} \label{remark:1-cocycle}
 Let  $T^t=T+\sum_{i=1}^{\infty}T^it^i$ be any formal deformation of    $T=\{T_\alpha\}_{\alpha\in\Omega}$. Then the infinitesimal $T^1$ is a 1-cocycle in the
cochain complex of $T$.  Moreover, the corresponding cohomology class depends
only on the equivalence class of the deformation $T^t$.

An $\Gamma$-twisted    Rota-Baxter family $T=\{T_\alpha:V\rightarrow L\}_{\alpha\in\Omega}$  is said to be rigid if any
formal deformation $T^t=T+\sum_{i=1}^{\infty}T^it^i$  of $T$ is equivalent to the deformation $\overline{T}^t=T$.
 \end{remark}

\begin{theorem}
Let   $T=\{T_\alpha:V\rightarrow L\}_{\alpha\in\Omega}$ be an $\Gamma$-twisted    Rota-Baxter family.
  If $H^{1 }_{\mathrm{\Omega}-\mathrm{LY}}(V,L)=0$, then $T$ is
rigid.
\end{theorem}

\begin{proof}
Let $T^t=T+\sum_{i=1}^{\infty}T^it^i$ be any formal deformation of    $T=\{T_\alpha\}_{\alpha\in\Omega}$. Then by Remark   \ref{remark:1-cocycle},
the  infinitesimal term $T^1$
is a 1-cocycle in the cohomology of $T$, i.e., $T^1\in Z^{1 }_{\mathrm{\Omega}-\mathrm{LY}}(V,L)$. Thus, by $H^{1 }_{\mathrm{\Omega}-\mathrm{LY}}(V,L)=0$, there exists $(a,b)\in C^{0 }_{\mathrm{\Omega}-\mathrm{LY}}(V,L)$ such that
 $T^1=\partial_\Omega(a,b)$. We take
 \begin{align*}
\zeta^t=\big\{\zeta^t_\alpha:= \mathrm{id}_V+t\big(D(a,b)(\cdot)+\Gamma_2(a,b,T_\alpha(\cdot))\big) \big\}_{\alpha\in \Omega},~  \ ~ \eta^t:=\mathrm{id}_L+t\{a,b,\cdot\}
\end{align*}
and define $\overline{T}^t= \eta^t\circ T^t\circ (\zeta^t)^{-1}$.  Then $\overline{T}^t$ is equivalent to the deformation ${T}^t$. Moreover, for all $u\in V, $  we get
  \begin{align*}
 \overline{T}^t(u)=&\Big\{(\mathrm{id}_L+t\{a,b,\cdot\})\circ(T_\alpha+tT_\alpha^1+\text{power of} ~t^{\geq 2})\big(u-t\big(D(a,b)u+\Gamma_2(a,b,T_\alpha u))\big)+\text{power of} ~t^{\geq 2}\Big\}_{\alpha\in \Omega}\\
 =&\Big\{(\mathrm{id}_L+t\{a,b,\cdot\})(T_\alpha u+t\big(T_\alpha^1u-T_\alpha \big(D(a,b)u+\Gamma_2(a,b,T_\alpha u))\big)+\text{power of} ~t^{\geq 2}\Big\}_{\alpha\in\Omega}\\
  =&\Big\{T_\alpha u+t\big(T_\alpha^1u+\{a,b,T_\alpha u\}-T_\alpha \big(D(a,b)u+\Gamma_2(a,b,T_\alpha u))\big)+\text{power of} ~t^{\geq 2}\Big\}_{\alpha\in \Omega}\\
  =&\Big\{T_\alpha u+t^2 \bar{{T}}^2_\alpha u+\cdots\Big\}_{\alpha\in \Omega},
 \end{align*}
where the last equation follows by $T^1=\partial_\Omega(a,b)$. This indicates that the coefficient
of $t$ in   $\overline{T}^t$ is trivial. By repeating this process, we get that $T^t$ is
equivalent to $T$,  indicating that $T$ is indeed rigid.
\end{proof}

{\bf Declaration of interests}

  The author declares no conflict of interest.

 {\bf Data availability}

No data was used for the research described in the article.

{{\bf Acknowledgments}

The paper is  supported by the Science and Technology Program of Guizhou Province (Grant No. QKHJC QN[2025]362).

\end{document}